\documentclass{article}
\usepackage{amsmath,bm,colortbl,graphicx,mathrsfs,afterpage,amssymb,url,cite,amsthm}
\RequirePackage[colorlinks,citecolor=blue,urlcolor=blue]{hyperref}
\pdfoutput=1
\def\P{{\rm P}} 
\def\E{{\rm E}} 
\allowdisplaybreaks

\begin{document}
\title{Parrondo games with spatial dependence, III}
\author{S. N. Ethier\thanks{Department of Mathematics, University of Utah, 155 South 1400 East, Salt Lake City, UT 84112, USA. e-mail: ethier@math.utah.edu.  Partially supported by a grant from the Simons Foundation (209632).}\ \ and 
Jiyeon Lee\thanks{Department of Statistics, Yeungnam University, 214-1 Daedong, Kyeongsan, Kyeongbuk 712-749, South Korea. e-mail: leejy@yu.ac.kr.  Supported by a 2012 Yeungnam University Research Grant.}}
\date{}
\maketitle

\begin{abstract}
We study Toral's Parrondo games with $N$ players and one-dimensional spatial dependence as modified by Xie et al.  Specifically, we use computer graphics to sketch the Parrondo and anti-Parrondo regions for $3\le N\le 9$.  Our work was motivated by a recent paper of Li et al., who applied a state space reduction method to this model, reducing the number of states from $2^N$ to $N+1$.  We show that their reduced Markov chains are inconsistent with the model of Xie et al.\medskip\par

\end{abstract}

\section{Introduction}\label{intro}

The Parrondo effect refers to a reversal in direction of some system parameter when two similar dynamics are combined.  It was first described by J. M. R. Parrondo in 1996 in the context of games of chance:  He showed that there exist two losing games that can be combined to win.  The games were originally intended as a pedagogical model of the flashing Brownian ratchet.  Early work focussed on capital-dependent (Harmer and Abbott \cite{HA99}) and history-dependent (Parrondo, Harmer, and Abbott \cite{PHA00}) games for a single player.  Multi-player games were introduced by Toral \cite{T01,T02}, including games with spatial dependence and games with redistribution of wealth.  

Toral's \cite{T01} Parrondo games with one-dimensional spatial dependence rely on an integer $N\ge3$ and three probability parameters, $p_0$, $p_1$, and $p_2$.  There are $N$ players labeled from 1 to $N$ and arranged in a circle in clockwise order.  At each turn, one player is chosen at random to play.  In game $B$, he tosses a $p_m$-coin if $m$ of his two nearest neighbors are winners ($m=0,1,2$).  A player's status as winner or loser depends on the result of his most recent game.  The player wins one unit with heads and loses one unit with tails.  The game can be initialized arbitrarily.  Game $A$ is the special case of game $B$ in which $p_0=p_1=p_2=1/2$.  For $0<\gamma<1$, game $C$, often denoted by $\gamma A+(1-\gamma)B$, is a random mixture of games $A$ and $B$.  (At each turn a coin with $\P(\text{heads})=\gamma$ is tossed, and game $A$ is played if heads appears, game $B$ if tails.)  We let $\mu_A$, $\mu_B$, and $\mu_C$ be the equilibrium mean profits per turn (to the ensemble of $N$ players) in games $A$, $B$, and $C$.  Of course, $\mu_A=0$.  We say that the \textit{Parrondo effect} occurs if $\mu_B\le0$ and $\mu_C>0$ (two fair or losing games combine to win) and the \textit{anti-Parrondo effect} occurs if $\mu_B\ge0$ and $\mu_C<0$ (two fair or winning games combine to lose).

Mihailovi\'c and Rajkovi\'c \cite{MR03} modeled the three games by introducing Markov chains in a state space with $2^N$ states.  They were able to study the games analytically for $N\le12$.  Ethier and Lee \cite{EL12a} used a state space reduction method that allowed them to study the games analytically for $N\le19$, an improvement that led to the conjecture that the mean profits $\mu_B$ and $\mu_C$ converge as $N\to\infty$.  (This was subsequently proved in \cite{EL13a} under certain conditions.)  Xie et al.\ \cite{X11} introduced a more spatially dependent version of game $A$, which we will refer to as game $A'$.  This game amounts to a loss of one unit by a randomly chosen player together with a win of one unit by a randomly chosen nearest neighbor of that player.  (This differs from the description in \cite{X11} but it is probabilistically equivalent.)  Li et al.\ \cite{L14} used a state space reduction method to study games $A'$, $B$, and $C':=\gamma A'+(1-\gamma)B$ analytically.  However, their method results in Markov chains that are not directly related to the model of Xie et al.

Let us explain what we mean by state space reduction, and what is required for it to be successful.

In general, consider an equivalence relation $\sim$ on a finite set $E$.  By definition, $\sim$ is \textit{reflexive} ($x\sim x$), \textit{symmetric} ($x\sim y$ implies $y\sim x$), and \textit{transitive} ($x\sim y$ and $y\sim z$ imply $x\sim z$).  It is well known that an equivalence relation partitions the set $E$ into \textit{equivalence classes}.  The set of all equivalence classes, called the \textit{quotient set}, will be denoted by $\bar E$.  Let us write $[x]:=\{y\in E: y\sim x\}$ for the equivalence class containing $x$.  Then $\bar E=\{[x]:x\in E\}$.

Now suppose $X_0,X_1,X_2,\ldots$ is a (time-homogeneous) Markov chain in $E$ with transition matrix $\bm P$.  In particular, $P(x,y)=\P(X_{t+1}=y\mid X_t=x)$ for all $x,y\in E$ and $t=0,1,2,\ldots$.  Under what conditions on $\bm P$ is $[X_0],[X_1],[X_2],\ldots$ a Markov chain in the ``reduced'' state space $\bar E$?  A sufficient condition, apparently due to Kemeny and Snell \cite[p.\ 124]{KS76}, is that $\bm P$ be \textit{lumpable} with respect to $\sim$.  By definition, this means that, for all $x,x',y\in E$,
\begin{equation}\label{lumpability}
x\sim x'\quad\text{implies}\quad \sum_{y'\in[y]}P(x,y')=\sum_{y'\in[y]}P(x',y').
\end{equation}
Moreover, if \eqref{lumpability} holds, then the Markov chain $[X_0],[X_1],[X_2],\ldots$ in $\bar E$ has transition matrix $\bar{\bm P}$ given by
\begin{equation}\label{Pbar}
\bar P([x],[y]):=\sum_{y'\in[y]}P(x,y').
\end{equation}
Notice that \eqref{lumpability} ensures that \eqref{Pbar} is well defined.

For Parrondo games with one-dimensional spatial dependence, the state space, assuming $N\ge3$ players, is
$$
E:=\{\bm x=(x_1,x_2,\ldots,x_N): x_i\in\{0,1\}{\rm\ for\ }i=1,\ldots,N\}=\{0,1\}^N,
$$
which has $2^N$ states.  A state $\bm x\in E$ describes the status of each of the $N$ players, 0 for losers and 1 for winners.  We can also think of $E$ as the set of $N$-bit binary representations of the integers $0,1,\ldots,2^N-1$, thereby giving a natural ordering to the vectors in $E$.

Ethier and Lee \cite{EL12a} used the following equivalence relation on $E$:  $\bm x\sim\bm y$ if and only if $\bm y=\bm x_\sigma:=(x_{\sigma(1)},\ldots,x_{\sigma(N)})$ for a permutation $\sigma$ of $(1,2,\ldots,N)$ belonging to the dihedral group $G$ of order $2N$ generated by the rotations and reflections of the players.  They verified the lumpability condition, with the result that the size of the state space was reduced by a factor of nearly $2N$ for large $N$.  It should be noted that a sufficient condition for the lumpability condition in this setting is
\begin{equation}\label{lump-alt}
P(\bm x,\bm y)=P(\bm x_\sigma,\bm y_\sigma)\quad\text{for all $\sigma\in G$}
\end{equation}
or for all $\sigma$ in a subset of $G$ that generates $G$.

Li et al.\ \cite{L14} reduced the state space much further by effectively using the following equivalence relation on $E$:  $\bm x\sim\bm y$ if and only if $\bm y=\bm x_\sigma$ for some permutation $\sigma$ of $(1,2,\ldots,N)$.  Equivalently, $\bm x\sim\bm y$ if and only if $\bm x$ and $\bm y$ have the same number of 1s.  However, here the lumpability condition fails.  (This is intuitively clear because the equivalence relation does not respect the spatial structure of the players.)  Thus, the Markov chains in $\bar E$ constructed in \cite{L14}, which are based on 
\begin{equation}\label{Pbar-ave}
\bar P([x],[y]):=\frac{1}{|[x]|}\sum_{x'\in[x]}\bigg[\sum_{y'\in[y]}P(x',y')\bigg]
\end{equation}
rather than \eqref{Pbar}, are not directly related to the model of Xie et al. \cite{X11}.  To put it another way, there is no theoretical basis for such Markov chains and the results derived from them.  

Our aim in this paper is not simply to point out a weakness in earlier work, but to provide accurate results as well for the model of Toral \cite{T01} as modified by Xie et al.\ \cite{X11}.  

Actually, as originally formulated, the model had four probability parameters $p_0$, $p_1$, $p_2$, and $p_3$ corresponding to the various configurations of the two nearest neighbors, namely $(0,0)$, $(0,1)$, $(1,0)$, and $(1,1)$.  However, we assume throughout that $p_1=p_2$ (and we relabel the remaining parameters $p_0,p_1,p_3$ as $p_0,p_1,p_2$).  There are three reasons for doing this.  First, the Parrondo region is a subset of the parameter space, and a three-dimensional region is easier to visualize than a four-dimensional one.  Second, the equivalence relation of Ethier and Lee \cite{EL12a} described above requires this assumption (otherwise the dihedral group of order $2N$ would have to be replaced by a cyclic group of order $N$ and the state space reduction would be less effective).  Third, all previous computational works on Parrondo games with one-dimensional spatial dependence \cite{T01,MR03,X11,EL12a,EL12b,L14} have made this assumption.

The smallest number of players for which there is a distinction between the equivalence relation of \cite{EL12a} and that of \cite{L14} is $N=4$, so we consider that case first in Section \ref{N=4}.  We treat the  general case by the method of \cite{EL12a} in Section \ref{EL_method} and by the method of \cite{L14} in Section \ref{Li_method}.  In Section \ref{Parrondo} we use computer graphics to sketch the Parrondo and anti-Parrondo regions for $3\le N\le 9$ as well as their approximations based on the methods of Li et al.\ \cite{L14}.  Finally, Section \ref{Conclusions} summarizes the main conclusions.

\section{The case $N=4$}\label{N=4}

To make this less abstract, let us consider separately the case of four players, that is, $N=4$.  We begin with game $B$.  The Markov chain in $E$ for game $B$ depends on three parameters, $p_0,p_1,p_2\in[0,1]$, so that the parameter space is the unit cube $[0,1]^3$.  There are 16 states (namely, the four-bit binary representations of the integers 0--15) and the transition matrix has the form
$$
\arraycolsep=1.1mm
\bm P_B:=\frac{1}{4}\left(
\begin{array}{cccccccccccccccc}
d_0 & p_0 & p_0 & 0 & p_0 & 0 & 0 & 0 & p_0 & 0 & 0 & 0 & 0 & 0 & 0 & 0 \\ 
q_0 & d_1 & 0 & p_1 & 0 & p_0 & 0 & 0 & 0 & p_1 & 0 & 0 & 0 & 0 & 0 & 0 \\ 
q_0 & 0 & d_2 & p_1 & 0 & 0 & p_1 & 0 & 0 & 0 & p_0 & 0 & 0 & 0 & 0 & 0 \\ 
0 & q_1 & q_1 & d_3 & 0 & 0 & 0 & p_1 & 0 & 0 & 0 & p_1 & 0 & 0 & 0 & 0 \\ 
q_0 & 0 & 0 & 0 & d_4 & p_0 & p_1 & 0 & 0 & 0 & 0 & 0 & p_1 & 0 & 0 & 0 \\ 
0 & q_0 & 0 & 0 & q_0 & d_5 & 0 & p_2 & 0 & 0 & 0 & 0 & 0 & p_2 & 0 & 0 \\ 
0 & 0 & q_1 & 0 & q_1 & 0 & d_6 & p_1 & 0 & 0 & 0 & 0 & 0 & 0 & p_1 & 0 \\ 
0 & 0 & 0 & q_1 & 0 & q_2 & q_1 & d_7 & 0 & 0 & 0 & 0 & 0 & 0 & 0 & p_2 \\ 
q_0 & 0 & 0 & 0 & 0 & 0 & 0 & 0 & d_8 & p_1 & p_0 & 0 & p_1 & 0 & 0 & 0 \\ 
0 & q_1 & 0 & 0 & 0 & 0 & 0 & 0 & q_1 & d_9 & 0 & p_1 & 0 & p_1 & 0 & 0 \\ 
0 & 0 & q_0 & 0 & 0 & 0 & 0 & 0 & q_0 & 0 & d_{10} & p_2 & 0 & 0 & p_2 & 0 \\ 
0 & 0 & 0 & q_1 & 0 & 0 & 0 & 0 & 0 & q_1 & q_2 & d_{11} & 0 & 0 & 0 & p_2 \\ 
0 & 0 & 0 & 0 & q_1 & 0 & 0 & 0 & q_1 & 0 & 0 & 0 & d_{12} & p_1 & p_1 & 0 \\ 
0 & 0 & 0 & 0 & 0 & q_2 & 0 & 0 & 0 & q_1 & 0 & 0 & q_1 & d_{13} & 0 & p_2 \\ 
0 & 0 & 0 & 0 & 0 & 0 & q_1 & 0 & 0 & 0 & q_2 & 0 & q_1 & 0 & d_{14} & p_2 \\ 
0 & 0 & 0 & 0 & 0 & 0 & 0 & q_2 & 0 & 0 & 0 & q_2 & 0 & q_2 & q_2 & d_{15}\end{array}\right),
$$
where the diagonal entries are chosen to make the row sums equal to 1:
\begin{eqnarray*}
d_0&:=&4 q_0,\\
d_1&=&d_2=d_4=d_8:=p_0 + q_0 + 2q_1,\\
d_3&=&d_6=d_9=d_{12}:=2 (p_1 + q_1),\\
d_5&=&d_{10}:=2 (p_0 + q_2),\\
d_7&=&d_{11}=d_{13}=d_{14}:=2 p_1 + p_2 + q_2,\\
d_{15}&:=&4 p_2,
\end{eqnarray*}
and $q_m:=1-p_m$ for $m=0,1,2$.  This is consistent with Eq.\ (12) of Xie et al.\ \cite{X11}.

For the equivalence relation of \cite{EL12a} mentioned above, there are six equivalence classes, namely
$\{0000\}$, $\{0001,0010,0100,1000\}$, $\{0011,0110,1001,1100\}$, $\{0101,1010\}$, $\{0111,1011,1101,1110\}$, and $\{1111\}$.
The lumpability condition is easily verified.  For example, denoting the states by their decimal representations (0--15),  the equivalence classes are $\{0\}$, $\{1,2,4,8\}$, $\{3,6,9,12\}$, $\{5,10\}$, $\{7,11,13,14\}$, and $\{15\}$, and the probabilities of a transition from a state in the equivalence class $\{3,6,9,12\}$ to the equivalence class itself are
\begin{eqnarray*}
\sum_{j=3,6,9,12}P_B(3,j)&=&(p_1 + q_1)/2,\quad
\sum_{j=3,6,9,12}P_B(6,j)=(p_1 + q_1)/2,\\
\sum_{j=3,6,9,12}P_B(9,j)&=&(p_1 + q_1)/2,\quad
\sum_{j=3,6,9,12}P_B(12,j)=(p_1 + q_1)/2,
\end{eqnarray*}
which are equal.  Moreover, the reduced transition matrix (with rows and columns labeled by the equivalence classes and ordered as just indicated) is 
\begin{equation}\label{pBbar-EL}
\setlength{\arraycolsep}{1.1mm}
\bar{\bm P}_B=\frac{1}{4}\left(\begin{array}{cccccc}
4 q_0 & 4 p_0 & 0 & 0 & 0 & 0 \\ 
q_0 & p_0 + q_0 + 2q_1 & 2p_1 & p_0 & 0 & 0 \\ 
0 & 2q_1 & 2 (p_1 + q_1) & 0 & 2p_1 & 0 \\ 
0 & 2 q_0 & 0 & 2 (p_0 + q_2) & 2 p_2 & 0 \\ 
0 & 0 & 2 q_1 & q_2 & 2 p_1 + p_2 + q_2 & p_2 \\ 
0 & 0 & 0 & 0 & 4 q_2 & 4 p_2
\end{array}\right).
\end{equation}
Both matrices, $\bm P_B$ and $\bar{\bm P}_B$, can be simplified (using $p_m + q_m =1$ for $m=0,1,2$), but there is a good reason to leave them in this form for now, as we will see. 

For the equivalence relation of \cite{L14} mentioned above, there are five equivalence classes, namely
$\{0000\}$, $\{0001,0010,0100,1000\}$, $\{0011,0101,0110,1001,1010,\break 1100\}$, $\{0111,1011,1101,1110\}$, and $\{1111\}$.
The lumpability condition fails.  For example, denoting the states by their decimal representations (0--15), the probabilities of a transition from a state in the equivalence class $\{3,5,6,9,10,12\}$ to the equivalence class itself are  
\begin{eqnarray*}
\sum_{j=3,5,6,9,10,12}P_B(3,j)&=&(p_1 + q_1)/2,\quad
\sum_{j=3,5,6,9,10,12}P_B(5,j)=(p_0 + q_2)/2,\\
\sum_{j=3,5,6,9,10,12}P_B(6,j)&=&(p_1 + q_1)/2,\quad
\sum_{j=3,5,6,9,10,12}P_B(9,j)=(p_1 + q_1)/2,\\
\sum_{j=3,5,6,9,10,12}P_B(10,j)&=&(p_0 + q_2)/2,\quad
\sum_{j=3,5,6,9,10,12}P_B(12,j)=(p_1 + q_1)/2,
\end{eqnarray*}
which are not, in general, equal.  The idea of Li et al. \cite{L14} is to use the average of these probabilities, namely $(p_0 + 2 p_1 + 2 q_1 + q_2)/6$, in $\bar{\bm P}_B$ (cf.\ \eqref{Pbar-ave}):
\begin{equation}\label{pBbar-Li}
\arraycolsep=0.2mm
\bar{\bm P}_B=\frac{1}{12}\left(
\begin{array}{ccccc}
12 q_0 & 12 p_0 & 0 & 0 & 0 \\ 
3 q_0 & 3 (p_0 + q_0 + 2 q_1) & 3 (p_0 + 2 p_1) & 0 & 0 \\ 
0 & 2 (q_0 + 2 q_1) & 2 (p_0 + 2 p_1 + 2 q_1 + q_2) & 2 (2 p_1 + p_2) & 0 \\ 
0 & 0 & 3 (2 q_1 + q_2) & 3 (2 p_1 + p_2 + q_2) & 3 p_2 \\ 
0 & 0 & 0 & 12 q_2 & 12 p_2 \end{array}\right),
\end{equation}
which is consistent with Eq.\ (3) of Li et al.\ \cite{L14}.

Next we need some notation.  We denote by $\dot{\bm P}_B$ the matrix $\bm P_B$ with $q_m$ replaced by $-q_m$ for $m=0,1,2$.  Similarly, we denote by $\dot{\bar{\bm P}}_B$ the matrix $\bar{\bm P}_B$ with $q_m$ replaced by $-q_m$ for $m=0,1,2$.  Let $\bm\pi_B$ be the unique stationary distribution of $\bm P_B$, and let $\bar{\bm\pi}_B$ be the unique stationary distribution of $\bar{\bm P}_B$.  (We assume that they exist, which rules out certain boundary cases such as $p_0=0$, $p_1$ arbitrary, and $p_2=1$.)  Then, as shown in \cite{EL12a} using the lumpability property (see \cite{EL12b} for more detail), the mean profit per turn in game $B$ is given by
\begin{equation}\label{ruleB}
\mu_B=\bm\pi_B\dot{\bm P}_B\bm1=\bar{\bm\pi}_B\dot{\bar{\bm P}}_B{\bm1},
\end{equation}
where $\bm 1$ is the column vector of 1s of the appropriate dimension. This equation shows that, as far as $\mu_B$ is concerned, nothing is lost in reducing the state space from $E$ to $\bar E$, provided the lumpability condition is satisfied (as in the case of \eqref{pBbar-EL}).  However, when this condition is not satisfied (as in the case of \eqref{pBbar-Li}), the relationship between $\bm P_B$ and $\bar{\bm P}_B$ is less clear, and in particular the second equation in \eqref{ruleB} fails.

Straightforward computations yield, for the equivalence relation of \cite{EL12a},
\begin{eqnarray}\label{muB-exact}
\mu_B&=&[-(3 - 2p_2 - 3p_0^2 + 2p_0 p_2 - p_2^2 + 2p_0^2 p_2 - 2p_0 p_2^2) \nonumber\\
 &&\quad{}+ 4 (1 + p_0) (1 - p_0 + p_2) (1 - p_2) p_1 - 
 2 (1 - p_0 + p_2) (1 - p_0 - p_2) p_1^2]\nonumber\\
 &&\;/[3 + 6 p_0 - 2 p_2 - 3 p_0^2 - 2 p_0 p_2 - p_2^2 + 12 p_0^2 p_2 - 4 p_0 p_2^2 - 8 p_0^2 p_2^2\nonumber\\
 &&\quad{} - 4 (1-p_0 + p_2 + 2 p_0^2 + 2 p_0 p_2) (1-p_2) p_1\nonumber\\
 &&\quad{} + 2 (1 + 4 p_0 - p_0^2 - 2 p_0 p_2 - p_2^2) p_1^2].
\end{eqnarray}
On the other hand, if we use the equivalence relation of \cite{L14}, we can compute the right-hand side of \eqref{ruleB}, but it yields only an approximation of $\mu_B$, which we denote by $\hat\mu_B$:
\begin{eqnarray}\label{muB-approx}
\hat\mu_B&=&[-(9 + 6 p_0 - 12 p_2 - 3 p_0^2 - 8 p_0 p_2 + 3 p_2^2 + 2 p_0^2 p_2 + 
 2 p_0 p_2^2) \nonumber\\
 &&\quad{}+ 2 (6 + 5 p_0 - 7 p_2 + p_0^2 - 4 p_0 p_2 + p_2^2) p_1 - 
 4 (1 - p_0 - p_2) p_1^2]\nonumber\\
 &&\;/[9 + 24 p_0 - 12 p_2 + 9 p_0^2 - 32 p_0 p_2 + 3 p_2^2 - 8 p_0^2 p_2 + 
 8 p_0 p_2^2\nonumber\\
&&\quad{}- 2 (6 - p_0 - 7 p_2 - p_0^2 + p_2^2) p_1 + 4 (1 + p_0 - p_2) p_1^2].
\end{eqnarray}
The fair surface $\mu_B=0$ based on \eqref{muB-exact} and its approximation $\hat\mu_B=0$ based on \eqref{muB-approx} are graphed in Figure \ref{ELvLi,B4,C4}(a).  The difference is significant.

\begin{figure}[ht]
\centering
\includegraphics[width = 2in]{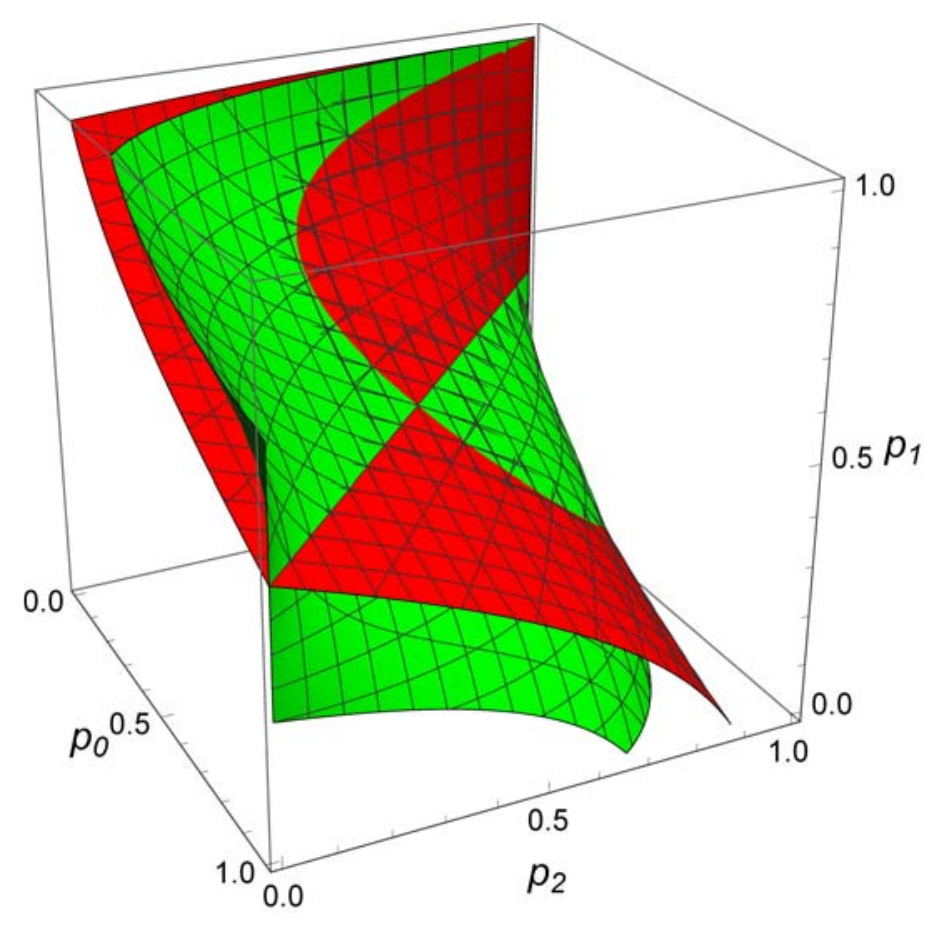}\qquad
\includegraphics[width = 2in]{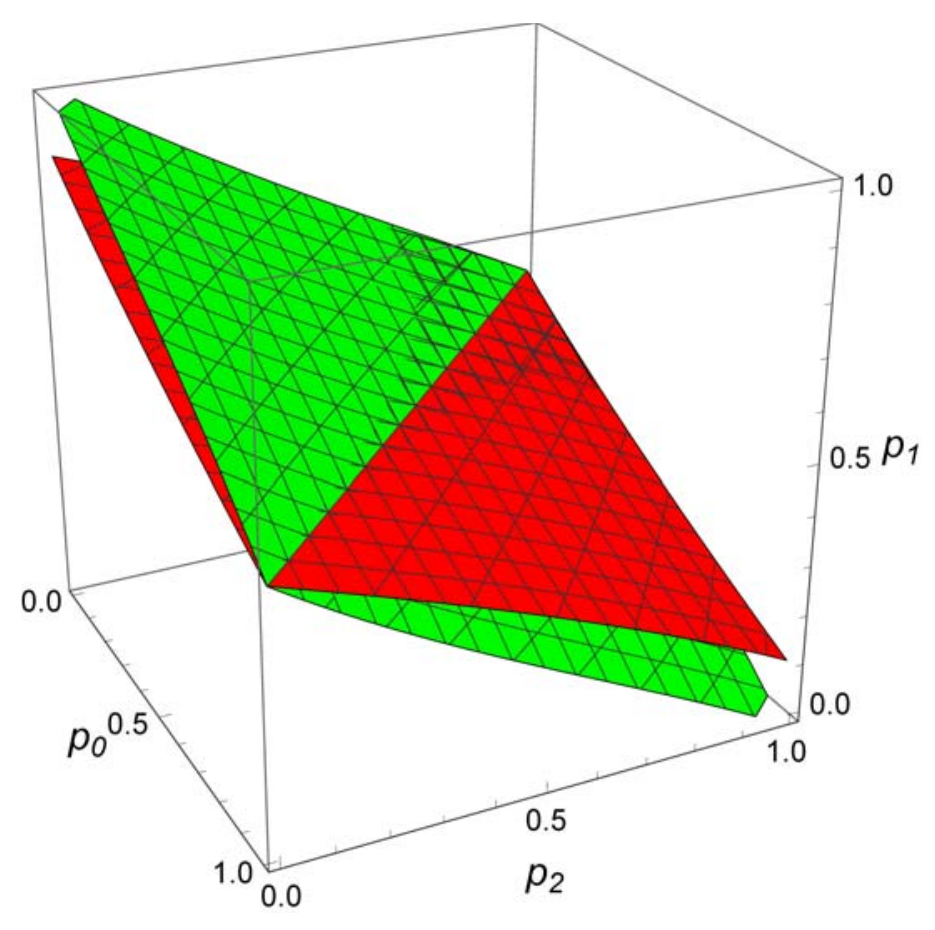}\\
\hglue-1cm (a)\hskip2.1in(b)
\caption{\label{ELvLi,B4,C4}(a) The fair surface $\mu_B=0$ for game $B$ with $N=4$ is displayed in green, while its approximation $\hat\mu_B=0$ is displayed in red.  (b) The fair surface $\mu_{C'}=0$ for game $C'$ with $N=4$ and $\gamma=1/2$ is displayed in green, while its approximation $\hat\mu_{C'}=0$ is displayed in red.}
\end{figure}

We turn next to game $A'$.  Again there are 16 states (namely, the 4-bit binary representations of the integers 0--15) and the transition matrix has the form
$$
\arraycolsep=2mm
\bm P_{A'}:=\frac{1}{8}\left(
\begin{array}{cccccccccccccccc}
0 & 2 & 2 & 0 & 2 & 0 & 0 & 0 & 2 & 0 & 0 & 0 & 0 & 0 & 0 & 0 \\ 
0 & 2 & 1 & 1 & 0 & 2 & 0 & 0 & 1 & 1 & 0 & 0 & 0 & 0 & 0 & 0 \\ 
0 & 1 & 2 & 1 & 1 & 0 & 1 & 0 & 0 & 0 & 2 & 0 & 0 & 0 & 0 & 0 \\ 
0 & 1 & 1 & 2 & 0 & 1 & 0 & 1 & 0 & 0 & 1 & 1 & 0 & 0 & 0 & 0 \\ 
0 & 0 & 1 & 0 & 2 & 2 & 1 & 0 & 1 & 0 & 0 & 0 & 1 & 0 & 0 & 0 \\ 
0 & 0 & 0 & 1 & 0 & 4 & 1 & 0 & 0 & 1 & 0 & 0 & 1 & 0 & 0 & 0 \\ 
0 & 0 & 1 & 0 & 1 & 1 & 2 & 1 & 0 & 0 & 1 & 0 & 0 & 0 & 1 & 0 \\ 
0 & 0 & 0 & 1 & 0 & 2 & 1 & 2 & 0 & 0 & 0 & 1 & 0 & 0 & 1 & 0 \\ 
0 & 1 & 0 & 0 & 1 & 0 & 0 & 0 & 2 & 1 & 2 & 0 & 1 & 0 & 0 & 0 \\ 
0 & 1 & 0 & 0 & 0 & 1 & 0 & 0 & 1 & 2 & 1 & 1 & 0 & 1 & 0 & 0 \\ 
0 & 0 & 0 & 1 & 0 & 0 & 1 & 0 & 0 & 1 & 4 & 0 & 1 & 0 & 0 & 0 \\ 
0 & 0 & 0 & 1 & 0 & 0 & 0 & 1 & 0 & 1 & 2 & 2 & 0 & 1 & 0 & 0 \\ 
0 & 0 & 0 & 0 & 1 & 1 & 0 & 0 & 1 & 0 & 1 & 0 & 2 & 1 & 1 & 0 \\ 
0 & 0 & 0 & 0 & 0 & 2 & 0 & 0 & 0 & 1 & 0 & 1 & 1 & 2 & 1 & 0 \\ 
0 & 0 & 0 & 0 & 0 & 0 & 1 & 1 & 0 & 0 & 2 & 0 & 1 & 1 & 2 & 0 \\ 
0 & 0 & 0 & 0 & 0 & 0 & 0 & 2 & 0 & 0 & 0 & 2 & 0 & 2 & 2 & 0\end{array}\right).
$$
The corresponding matrix in Eq.\ (3) of Xie et al.\ \cite{X11} has errors at six of the 256 entries.  Labeling rows and columns by 0--15, the errors occur at entries $(5,3)$, $(5,12)$, $(10,2)$, $(10,12)$, $(15,7)$, and $(15,9)$.

Using the equivalence relation of \cite{EL12a} with six equivalence classes, the lumpability condition holds and the reduced transition matrix (with rows and columns labeled by the equivalence classes and ordered as previously indicated) is 
\begin{equation}\label{pAbar-EL}
\setlength{\arraycolsep}{2mm}
\bar{\bm P}_{A'}=\frac{1}{4}\left(
\begin{array}{cccccc}
0 & 4 & 0 & 0 & 0 & 0 \\
0 & 2 & 1 & 1 & 0 & 0 \\
0 & 1 & 1 & 1 & 1 & 0 \\
0 & 0 & 2 & 2 & 0 & 0 \\
0 & 0 & 1 & 1 & 2 & 0 \\
0 & 0 & 0 & 0 & 4 & 0 \\
\end{array}\right).
\end{equation}

Using the equivalence relation of \cite{L14} with five equivalence classes, the lumpability condition fails, but using \eqref{Pbar-ave} yields the transition matrix
\begin{equation}\label{pAbar-Li}
\setlength{\arraycolsep}{2mm}
\bar{\bm P}_{A'}=\frac{1}{6}\left(
\begin{array}{lllll}
0 & 6 & 0 & 0 & 0 \\ 
0 & 3 & 3 & 0 & 0 \\ 
0 & 1 & 4 & 1 & 0 \\ 
0 & 0 & 3 & 3 & 0 \\ 
0 & 0 & 0 & 6 & 0
\end{array}\right)
\end{equation}
as in Eq.\ (1) of Li et al.\ \cite{L14}.

The transition matrices corresponding to game $C':=\gamma A'+(1-\gamma)B$ are 
$$
\bm P_{C'}:=\gamma\bm P_{A'}+(1-\gamma)\bm P_B \quad\text{and}\quad 
\bar{\bm P}_{C'}:=\gamma\bar{\bm P}_{A'}+(1-\gamma)\bar{\bm P}_B.
$$
Let $\bm\pi_{C'}$ be the unique stationary distribution of $\bm P_{C'}$, and let $\bar{\bm\pi}_{C'}$ be the unique stationary distribution of $\bar{\bm P}_{C'}$.  Since every play of game $A'$ results in a profit of 0 to the ensemble of four players, $\dot{\bm P}_{A'}=\bm0$ and $\dot{\bar{\bm P}}_{A'}=\bm0$ and hence
$$
\dot{\bm P}_{C'}=(1-\gamma)\dot{\bm P}_B \quad\text{and}\quad
\dot{\bar{\bm P}}_{C'} =(1-\gamma)\dot{\bar{\bm P}}_B.
$$
This allows us to evaluate
\begin{equation}\label{ruleC'}
\mu_{C'}=\bm\pi_{C'}\dot{\bm P}_{C'}\bm1=\bar{\bm\pi}_{C'}\dot{\bar{\bm P}}_{C'}{\bm1}.
\end{equation}

Let us assume that $\gamma=1/2$.  Using the equivalence relation of \cite{EL12a}, we find from \eqref{ruleC'} that
\begin{eqnarray}\label{muC-exact}
\mu_{C'}&=&[-3 (105 - 35 p_0 - 65 p_2 - 22 p_0^2 + 8 p_0 p_2 + 2 p_2^2 + 6 p_0^2 p_2 - 
      2 p_0 p_2^2)\nonumber\\
&&\quad{}   + 6 (55 + 2 p_0 - 4 p_2 - 9 p_0^2 + 12 p_0 p_2 - 9 p_2^2 + 4 p_0^2 p_2 - 
      4 p_0 p_2^2) p_1\nonumber\\
&&\quad{}   - 12 (2 - p_0 + p_2) (1 - p_0 - p_2) p_1^2]\nonumber\\
&&\;{}   /[2 (315 + 175 p_0 - 125 p_2 - 
      22 p_0^2 - 10 p_0 p_2 - 12 p_2^2 + 48 p_0^2 p_2 - 32 p_0 p_2^2 \nonumber\\
&&\quad{}  - 16 p_0^2 p_2^2) - 4 (25 - 2 p_0 + 8 p_2 + 17 p_0^2 + 12 p_0 p_2 - 13 p_2^2 - 8 p_0^2 p_2 \nonumber\\
&&\quad{} - 8 p_0 p_2^2) p_1  + 8 (14 + 7 p_0 - 3 p_2 - p_0^2 - 2 p_0 p_2 - p_2^2) p_1^2],
\end{eqnarray}
whereas, using the equivalence relation of \cite{L14}, the right-hand side of \eqref{ruleC'} yields the approximation
\begin{eqnarray}\label{muC-approx}
\hat\mu_{C'}&=& [-3 (56 - p_0 - 45 p_2 - 7 p_0^2 - 8 p_0 p_2 + 7 p_2^2 + 2 p_0^2 p_2 + 
      2 p_0 p_2^2)\nonumber\\
&&\quad{} + 6 (33 + 14 p_0 - 16 p_2 + p_0^2 - 4 p_0 p_2 + p_2^2) p_1 - 12 (1 - p_0 - p_2) p_1^2]\nonumber\\
&&\;{}/[2 (168 + 137 p_0 - 107 p_2 + 19 p_0^2 - 
      80 p_0 p_2 + 13 p_2^2 - 8 p_0^2 p_2 + 8 p_0 p_2^2)\nonumber\\
&&\quad{}  - 4 (15 - 6 p_0 - 12 p_2 - p_0^2 + p_2^2) p_1 + 8 (3 + p_0 - p_2) p_1^2].
\end{eqnarray}
The fair surface $\mu_{C'}=0$ based on \eqref{muC-exact} and its approximation $\hat\mu_{C'}=0$ based on \eqref{muC-approx} are graphed in Figure \ref{ELvLi,B4,C4}(b) above.  The difference is again significant.

\begin{figure}[ht]
\centering
\includegraphics[width = 2in]{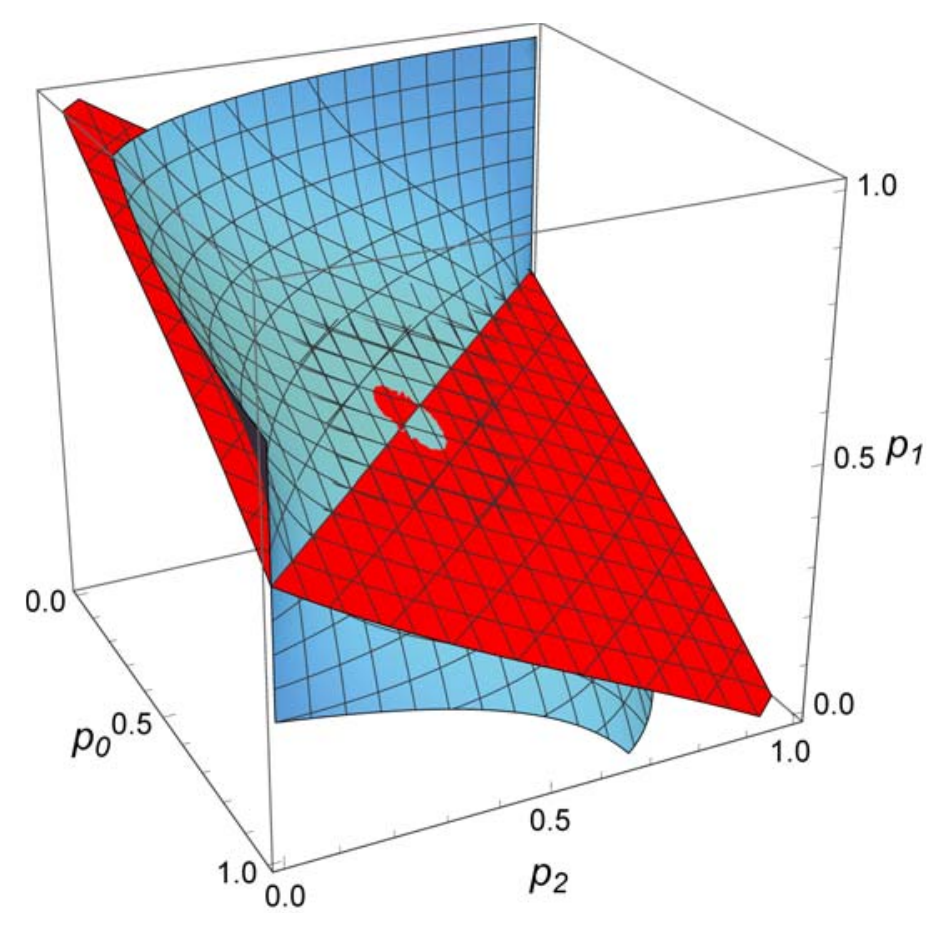}\quad
\includegraphics[width = 2in]{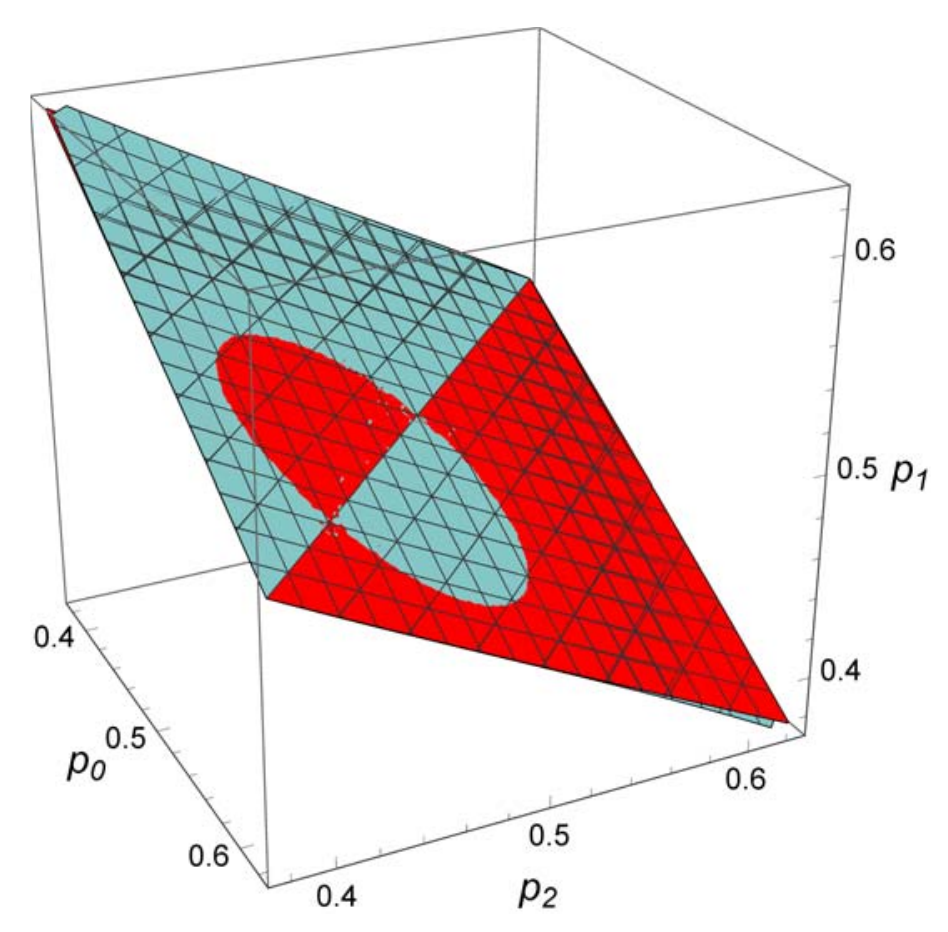}\\
\caption{\label{region_N=4}For $N=4$ and $\gamma=1/2$, the blue surface is the surface $\mu_B=0$, and the red surface is the surface $\mu_{C'}=0$, in the $(p_0,p_2,p_1)$ unit cube.  The Parrondo region is the region on or below the blue surface and above the red surface, while the anti-Parrondo region is the region on or above the blue surface and below the red surface.  The second figure is an enlargement of a portion of the first.}
\end{figure}

We can use the formulas found above to sketch the Parrondo and anti-Parrondo regions for games $A'$, $B$, and $C'$.  This is done in Figure \ref{region_N=4}.  We notice a curious feature of the regions that did not appear for games $A$, $B$, and $C$.   When projected onto the $(p_0,p_2)$ plane, there is a small apparently circular region in which the Parrondo and anti-Parrondo regions on opposite sides on the line $p_0+p_2=1$ are inverted.  However, this region is neither circular nor elliptical, but rather is determined by a polynomial of degree 8 in $p_0$ and $p_2$.

\section{Reduction by reflections and rotations}\label{EL_method}

The Markov chain formalized by Mihailovi\'c and Rajkovi\'c \cite{MR03} keeps track of the status (loser or winner, 0 or 1) of each of the $N\ge3$ players of game $B$, which was described in Section \ref{intro}. Its state space is the product space 
$$
E:=\{\bm x=(x_1,x_2,\ldots,x_N): x_i\in\{0,1\}{\rm\ for\ }i=1,\ldots,N\}=\{0,1\}^N
$$
with $2^N$ states.  Let $m_i(\bm x):=x_{i-1}+x_{i+1}$, or, in other words, let $m_i(\bm x)$ be the number of winners among the two nearest neighbors of player $i$ when the state is $\bm x$.  Of course $x_0:=x_N$ and $x_{N+1}:=x_1$ because of the circular arrangement of players.  Also, let $\bm x^i$ be the element of $E$ equal to $\bm x$ except at the $i$th component.  For example, $\bm x^1:=(1-x_1,x_2,x_3,\ldots,x_N)$.

The one-step transition matrix $\bm P_B$ for this Markov chain depends not only on $N$ but on three parameters, $p_0,p_1,p_2\in[0,1]$.  It has the form
\begin{equation*}
P_B(\bm x,\bm x^i):=\begin{cases}N^{-1}p_{m_i(\bm x)}&\text{if $x_i=0$,}\\N^{-1}q_{m_i(\bm x)}&\text{if $x_i=1$,}\end{cases}\qquad i=1,\ldots,N,\;\bm x\in E,
\end{equation*}
and 
\begin{equation*}
P_B(\bm x,\bm x):=N^{-1}\bigg(\sum_{i:x_i=0}q_{m_i(\bm x)}+\sum_{i:x_i=1}p_{m_i(\bm x)}\bigg),\qquad \bm x\in E,
\end{equation*}
where $q_m:=1-p_m$ for $m=0,1,2$ and empty sums are 0.  The Markov chain is irreducible and aperiodic if $0<p_m<1$ for $m=0,1,2$.  Under slightly weaker assumptions (see \cite{EL13a}), the Markov chain is ergodic, which suffices.  For example, if $p_0=1$ and $0<p_m<1$ for $m=1,2$, or if $0<p_m<1$ for $m=0,1$ and $p_2=0$, then ergodicity holds.

The Markov chain of Xie et al.\ \cite{X11} keeps track of the status (loser or winner, 0 or 1) of each of the $N\ge3$ players of game $A'$.  Recall that a player is chosen at random, and he then loses one unit while a randomly chosen nearest neighbor of that player wins one unit.  Its state space is again $E$.  For $i=1,2,\ldots,N$, let $\bm x^{i,1}$ be the element of $E$ whose $j$th component is equal to $x_j$ if $j\ne i,i+1$, 0 if $j=i$, and 1 if $j=i+1$; and let $\bm x^{i,-1}$ be the element of $E$ whose $j$th component is equal to $x_j$ if $j\ne i,i-1$, 0 if $j=i$, and 1 if $j=i-1$.  Of course $x_0:=x_N$ and $x_{N+1}:=x_1$.  For example, $\bm x^{1,-1}:=(0,x_2,x_3,\ldots,x_{N-1},1)$.  Starting from state $\bm x\in E$, player $i$ is chosen with probability $1/N$, and a transition to state $\bm x^{i,1}$ or state $\bm x^{i,-1}$ occurs, each with probability 1/2. 

The one-step transition matrix $\bm P_{A'}$ for this Markov chain in $E$ therefore has the form
\begin{equation*}
P_{A'}(\bm x,\bm y):=(2N)^{-1}\sum_{i=1}^N[\delta(\bm x^{i,-1},\bm y)+\delta(\bm x^{i,1},\bm y)],
\end{equation*}
where $\delta(\bm x,\bm y)$ is the Kronecker delta ($\delta(\bm x,\bm y):=1$ if $\bm x=\bm y$; $\delta(\bm x,\bm y):=0$ otherwise).  The Markov chain is not irreducible because states $\bm 0$ and $\bm 1$ cannot be reached, but it does have a unique stationary distribution because it is irreducible when restricted to $E-\{\bm0,\bm1\}$.

Using the equivalence relation of Ethier and Lee \cite{EL12a}, we can verify the lumpability condition \eqref{lump-alt} for the transition matrix $\bm P_B$ by checking that
\begin{eqnarray*}
P_B(\bm x_\sigma,(\bm x^i)_\sigma)&=&P_B(\bm x_\sigma,(\bm x_\sigma)^{\sigma^{-1}(i)})\\
&=&\begin{cases}N^{-1}p_{m_{\sigma^{-1}(i)}(\bm x_\sigma)}&\text{if $(x_\sigma)_{\sigma^{-1}(i)}=0$,}\\N^{-1}q_{m_{\sigma^{-1}(i)}(\bm x_\sigma)}&\text{if $(x_\sigma)_{\sigma^{-1}(i)}=1$,}\end{cases}\\
&=&\begin{cases}N^{-1}p_{m_i(\bm x)}&\text{if $x_i=0$,}\\N^{-1}q_{m_i(\bm x)}&\text{if $x_i=1$,}\end{cases}\\
&=&P_B(\bm x,\bm x^i),
\end{eqnarray*}
where the third equality uses
\begin{equation}\label{m-sigma}
m_i({\bm x}_\sigma)=m_{\sigma(i)}(\bm x).
\end{equation}
We can verify that \eqref{m-sigma} holds
for $(\sigma(1),\ldots,\sigma(N))=(2,3,\ldots,N,1)$ and for $(\sigma(1),\ldots,\sigma(N))=(N,N-1,\ldots,2,1)$, which suffices because these two permutations generate $G$.  

For $\bm P_{A'}$, we can verify the lumpability condition \eqref{lump-alt} by observing that, 
if $(\sigma(1),\ldots,\sigma(N))=(2,3,\ldots,N,1)$, then
\begin{eqnarray*}
P_{A'}(\bm x_\sigma,\bm y_\sigma)&=&(2N)^{-1}\sum_{i=1}^N[\delta((\bm x_\sigma)^{i,-1},\bm y_\sigma)+\delta((\bm x_\sigma)^{i,1},\bm y_\sigma)]\\
&=&(2N)^{-1}\sum_{i=1}^N[\delta((\bm x_\sigma)^{\sigma^{-1}(i),-1},\bm y_\sigma)+\delta((\bm x_\sigma)^{\sigma^{-1}(i),1},\bm y_\sigma)]\\
&=&(2N)^{-1}\sum_{i=1}^N[\delta((\bm x^{i,-1})_\sigma,\bm y_\sigma)+\delta((\bm x^{i,1})_\sigma,\bm y_\sigma)]\\
&=&(2N)^{-1}\sum_{i=1}^N[\delta(\bm x^{i,-1},\bm y)+\delta(\bm x^{i,1},\bm y)]\\
&=&P_{A'}(\bm x,\bm y)
\end{eqnarray*}
since $(\bm x_\sigma)^{\sigma^{-1}(i),s}=(\bm x^{i,s})_\sigma$.  If $(\sigma(1),\ldots,\sigma(N))=(N,N-1,\ldots,2,1)$, then the same sequence of identities holds, the only distinction being that $(\bm x_\sigma)^{\sigma^{-1}(i),s}=(\bm x^{i,-s})_\sigma$.
Since the lumpability condition is satisfied for $\bm P_{A'}$, it is also satisfied for $\bm P_{C'}$.

A fairly explicit formula for $\bar{\bm P}_B$ is given in \cite{EL12a}.  First, define the function $s:\bar{E}\mapsto\{0,1,\ldots,N\}$ by $s([\bm x]):=x_1+x_2+\cdots+x_N$; it counts the number of 1s in each element of an equivalence class.  Then
\begin{equation*}
\bar{P}_B([\bm x],[\bm y])=\begin{cases}N^{-1}\big(\sum_{i:x_i=0}q_{m_i(\bm x)}+\sum_{i:x_i=1}p_{m_i(\bm x)}\big)&\text{if $[\bm y]=[\bm x]$}\\
N^{-1}\sum_{i:x_i=1,\bm x^i\sim\bm y}q_{m_i(\bm x)}&\text{if $s([\bm y])=s([\bm x])-1$}\\
N^{-1}\sum_{i:x_i=0,\bm x^i\sim\bm y}p_{m_i(\bm x)}&\text{if $s([\bm y])=s([\bm x])+1$}\\
0&\text{otherwise}\end{cases}
\end{equation*}
for all $[\bm x],[\bm y]\in \bar{E}$.  This generalizes \eqref{pBbar-EL}.  
A formula for $\bar{\bm P}_{A'}$ is
\begin{equation*}
\bar{P}_{A'}([\bm x],[\bm y])=\begin{cases}(2N)^{-1}\sum_{i=1}^N\big(1_{[\bm y]}(\bm x^{i,-1})+1_{[\bm y]}(\bm x^{i,1})\big)&\text{if $|s([\bm y])-s([\bm x])|\le1$}\\
0&\text{otherwise}\end{cases}
\end{equation*}
for all $[\bm x],[\bm y]\in \bar{E}$.  This generalizes \eqref{pAbar-EL}.

As in Section \ref{N=4}, let $\bm\pi_B$ be the unique stationary distribution of $\bm P_B$, and let $\bar{\bm\pi}_B$ be the unique stationary distribution of $\bar{\bm P}_B$.  We denote by $\dot{\bm P}_B$ the matrix $\bm P_B$ with $q_m$ replaced by $-q_m$ for $m=0,1,2$.  Similarly, we denote by $\dot{\bar{\bm P}}_B$ the matrix $\bar{\bm P}_B$ with $q_m$ replaced by $-q_m$ for $m=0,1,2$.  Then the mean profit per turn in game $B$ is given by
\begin{equation*}
\mu_B=\bm\pi_B\dot{\bm P}_B\bm1=\bar{\bm\pi}_B\dot{\bar{\bm P}}_B{\bm1}.
\end{equation*}
The transition matrices corresponding to game $C':=\gamma A'+(1-\gamma)B$ are 
$$
\bm P_{C'}:=\gamma\bm P_{A'}+(1-\gamma)\bm P_B \quad\text{and}\quad 
\bar{\bm P}_{C'}:=\gamma\bar{\bm P}_{A'}+(1-\gamma)\bar{\bm P}_B.
$$
Let $\bm\pi_{C'}$ be the unique stationary distribution of $\bm P_{C'}$, and let $\bar{\bm\pi}_{C'}$ be the unique stationary distribution of $\bar{\bm P}_{C'}$.  Since every play of game $A'$ results in a profit of 0 to the ensemble of $N$ players, $\dot{\bm P}_{A'}=\bm0$ and $\dot{\bar{\bm P}}_{A'}=\bm0$ and hence
$$
\dot{\bm P}_{C'}=(1-\gamma)\dot{\bm P}_B \quad\text{and}\quad
\dot{\bar{\bm P}}_{C'} =(1-\gamma)\dot{\bar{\bm P}}_B.
$$
Then the mean profit per turn in game $C'$ is given by
\begin{equation}\label{ruleC1}
\mu_{C'}=\bm\pi_{C'}\dot{\bm P}_{C'}\bm1=\bar{\bm\pi}_{C'}\dot{\bar{\bm P}}_{C'}{\bm1}.
\end{equation}

There is a technical issue that must be addressed to fully justify our results.  The strong law of large numbers of Ethier and Lee \cite{EL09} does not apply directly because the payoffs are not completely specified by the one-step transitions of the Markov chain.  (For example, a transition from $\bm x$ to $\bm x$ could correspond to a payoff of $-1$ or $1$ in game $B$; $-1$, $0$, or $1$ in game $C'$.)  We addressed this issue for game $B$ in \cite{EL12a} by augmenting the state space.  A similar approach, but with a different augmentation, works for game $C'$.  We let $E^*:=E\times\{1,2,\ldots,N\}\times\{1,2\}$.  The state is $(\bm x,i,g)\in E^*$ if $\bm x$ describes the status of each player, $i$ is the label of the next player to play, and $g$ is the next  game to be played ($g=1$ for game $A$ and $g=2$ for game $B$).  The new one-step transition matrix $\bm P_{C'}^*$ has the form
\begin{eqnarray*}
P_{C'}^*((\bm x,i,2),(\bm x^i,j,h))&=&\begin{cases}N^{-1}\gamma(h)p_{m_i(\bm x)}&\text{if $x_i=0$},\\
N^{-1}\gamma(h)q_{m_i(\bm x)}&\text{if $x_i=1$},\end{cases}\\
P_{C'}^*((\bm x,i,2),(\bm x,j,h))&=&\begin{cases}N^{-1}\gamma(h)q_{m_i(\bm x)}&\text{if $x_i=0$},\\
N^{-1}\gamma(h)p_{m_i(\bm x)}&\text{if $x_i=1$},\end{cases}\\
P_{C'}^*((\bm x,i,1),(\bm y,j,h))&=&N^{-1}\gamma(h)[\delta(\bm x^{i,-1},\bm y)+\delta(\bm x^{i,1},\bm y)]/2,
\end{eqnarray*}
where $\gamma(1):=\gamma$ and $\gamma(2):=1-\gamma$.  For each transition there is an associated payoff, specifically
$$
w((\bm x,i,2),(\bm y,j,h)):=\begin{cases}1&\text{if ($\bm y=\bm x^i$, $x_i=0$) or ($\bm y=\bm x$, $x_i=1$),}\\
-1&\text{if ($\bm y=\bm x$, $x_i=0$) or ($\bm y=\bm x^i$, $x_i=1$),}\end{cases}
$$
and $w((\bm x,i,1),(\bm y,j,h))=0$.  This allows us to show that $\mu_{C'}=\bm\pi_{C'}^*\dot{\bm P}_{C'}^*\bm1$ using the SLLN and $\bm\pi_{C'}^*\dot{\bm P}_{C'}^*\bm1=\bm\pi_{C'}\dot{\bm P}_{C'}\bm1=\bar{\bm\pi}_{C'}\dot{\bar{\bm P}}_{C'}{\bm1}$, from which \eqref{ruleC1} follows.

\section{Reduction to the number of winners}\label{Li_method}

Assuming the equivalence relation in which two states in $E$ are equivalent if they have the same number of 1s, Li et al.\ \cite{L14} found some rather elegant formulas for the matrices $\bar{\bm P}_{A'}$ and $\bar{\bm P}_B$ (defined using \eqref{Pbar-ave}).  Both are tridiagonal matrices with rows and columns indexed by $\bar{E}:=\{0,1,2,\ldots,N\}$.

Specifically,
\begin{equation}\label{pA'(i,i+1)}
\bar{P}_{A'}(i,i+1)=\frac{\binom{N-2}{i}}{\binom{N}{i}},\quad i=0,1,\ldots,N-1,
\end{equation}
\begin{equation*}
\bar{P}_{A'}(i,i-1)=\frac{\binom{N-2}{i-2}}{\binom{N}{i}},\quad i=1,2,\ldots,N,
\end{equation*}
where $\binom{N-2}{N-1}=\binom{N-2}{-1}:=0$,
and
\begin{equation}\label{complement}
\bar{P}_{A'}(i,i)=1-\bar{P}_{A'}(i,i+1)-\bar{P}_{A'}(i,i-1),\quad i=1,2,\ldots,N-1,
\end{equation}
and $\bar{P}_{A'}(0,0)=\bar{P}_{A'}(N,N)=0$.  This generalizes \eqref{pAbar-Li}.

Moreover,
\begin{equation*}
\bar{P}_B(i,i+1)=\frac{\binom{N-3}{i}p_0+2\binom{N-3}{i-1}p_1+\binom{N-3}{i-2}p_2}{\binom{N}{i}},\quad i=0,1,\ldots,N-1,
\end{equation*}
\begin{equation*}
\bar{P}_B(i,i-1)=\frac{\binom{N-3}{i-1}q_0+2\binom{N-3}{i-2}q_1+\binom{N-3}{i-3}q_2}{\binom{N}{i}},\quad i=1,2,\ldots,N,
\end{equation*}
where $\binom{N-3}{N-1}=\binom{N-3}{N-2}=\binom{N-3}{-1}=\binom{N-3}{-2}:=0$
and $q_m:=1-p_m$ for $m=0,1,2$; finally,
\begin{eqnarray}\label{i to i}
\bar{P}_B(i,i)&=&\frac{\binom{N-3}{i}q_0+2\binom{N-3}{i-1}q_1+\binom{N-3}{i-2}q_2}{\binom{N}{i}}+\frac{\binom{N-3}{i-1}p_0+2\binom{N-3}{i-2}p_1+\binom{N-3}{i-3}p_2}{\binom{N}{i}},\nonumber\\
&&{}\qquad\qquad\qquad\qquad\qquad\qquad\qquad\qquad\qquad i=1,2,\ldots,N-1,\qquad
\end{eqnarray}
and $\bar{P}_B(0,0)=q_0$ and $\bar{P}_B(N,N)=p_2$.  (Defining $\bar{P}_B(i,i)$ directly instead of as a complementary probability as in \eqref{complement} simplifies the definition of $\dot{\bar{\bm P}}_B$.)  This generalizes \eqref{pBbar-Li}.

Since Li et al.\ \cite{L14} did not include proofs, we sketch a proof of \eqref{pA'(i,i+1)}:  Consider a uniformly distributed sequence of 0s and 1s of length $N$.  Then
\begin{eqnarray*}
&&\!\!\!\!\!\!\!\!\!\!\!\!\!\!\bar{P}_{A'}(i,i+1)\\
&=&(1/2)\E[\text{proportion of positions with a 0 preceded by a 0}\mid i\text{ 1s}]\\
&&\quad{}+(1/2)\E[\text{proportion of positions with a 0 followed by a 0}\mid i\text{ 1s}]\\
&=&(2N)^{-1}\E[\text{number of positions with a 0 preceded by a 0}\mid i\text{ 1s}]\\
&&\quad{}+(2N)^{-1}\E[\text{number of positions with a 0 followed by a 0}\mid i\text{ 1s}]\\
&=&\P(\text{0s in positions $N-1$ and $N$} \mid i\text{ 1s})\\
&=&\frac{\binom{N-2}{i}}{\binom{N}{i}},
\end{eqnarray*}
where the third equality uses an exchangeability argument.  The other formulas are proved analogously.

These equations determine $\bar{\bm P}_{C'}:=\gamma\bar{\bm P}_{A'}+(1-\gamma)\bar{\bm P}_B$.  As before, 
we denote by $\dot{\bar{\bm P}}_B$ the matrix $\bar{\bm P}_B$ with $q_m$ replaced by $-q_m$ for $m=0,1,2$.  Since every play of game $A'$ results in a profit of 0 to the ensemble of $N$ players, $\dot{\bar{\bm P}}_{A'}=\bm0$ and hence
$$
\dot{\bar{\bm P}}_{C'}=(1-\gamma)\dot{\bar{\bm P}}_B.
$$
Let $\bar{\bm\pi}_B$ be the unique stationary distribution of $\bar{\bm P}_B$ on $\bar E$, and
let $\bar{\bm\pi}_{C'}$ be the unique stationary distribution of $\bar{\bm P}_{C'}$ on $\bar E$.  This allows us to evaluate
\begin{equation}\label{ruleLi}
\hat\mu_B:=\bar{\bm\pi}_B\dot{\bar{\bm P}}_B{\bm1}\quad\text{and}\quad
\hat\mu_{C'}:=\bar{\bm\pi}_{C'}\dot{\bar{\bm P}}_{C'}{\bm1}.
\end{equation}
These are the approximations of Li et al.\ \cite{L14} to the exact $\mu_B$ and $\mu_{C'}$.  

As in Section \ref{EL_method}, the strong law of large numbers in \cite{EL09} does not apply directly, so the Markov chains in $\bar{E}$ described above must be augmented.  We explain the procedure for game $B$; the procedure for game $C'$ is similar.  We let $\bar{E}^*:=\bar{E}\times\{-1,1\}$.  The state is $(i,w)\in \bar{E}^*$ if $i$ is the number of 1s and $w$ is the profit obtained in the transition to state $i$.  In \eqref{i to i}, notice that $\bar{P}_B(i,i)$ is the sum of two fractions, which we denote by $\bar{P}_B'(i,i)$ and $\bar{P}_B''(i,i)$, respectively.  The new one-step transition matrix $\bar{\bm P}_B^*$ has the form
\begin{eqnarray*}
{\bar P}_B^*((i,w),(i+1,1))&=&\bar{P}_B(i,i+1),\\
{\bar P}_B^*((i,w),(i+1,-1))&=&0,\\
{\bar P}_B^*((i,w),(i-1,1))&=&0,\\
{\bar P}_B^*((i,w),(i-1,-1))&=&\bar{P}_B(i,i-1),\\
{\bar P}_B^*((i,w),(i,1))&=&\bar{P}_B''(i,i),\\
{\bar P}_B^*((i,w),(i,-1))&=&\bar{P}_B'(i,i).
\end{eqnarray*}
This allows us to show that $\hat\mu_B=\bar{\bm\pi}_B^*\dot{\bar{\bm P}}_B^*\bm1$ using the SLLN and $\bar{\bm\pi}_B^*\dot{\bar{\bm P}}_B^*\bm1=\bar{\bm\pi}_B\dot{\bar{\bm P}}_B\bm1$, from which the first equation in \eqref{ruleLi} follows.

\section{The Parrondo region}\label{Parrondo}

We can now compare the Parrondo regions for games $A'$, $B$, and $C'$ of Xie et al.\ \cite{X11} with those for games $A$, $B$, and $C$ of Toral \cite{T01}.  We can also compare the Parrondo regions for games $A'$, $B$, and $C'$ computed exactly with those computed using the approximation of Li et al.\ \cite{L14}.  Let us first compare them for certain parameter values.  Table \ref{Toral} assumes that $(p_0,p_1,p_2)=(1,4/25,7/10)$, as did Toral \cite{T01}.  We can compute mean profits for $3\le N\le 19$.  Notice that the Parrondo effect appears for games $A$, $B$, and $C$, as well as for games $A'$, $B$, and $C'$, when $N=5,6$ and $N\ge9$.  However, the approximate formulas give a very different conclusion: There is no effect for $N\le15$, an anti-Parrondo effect for $16\le N\le28$, and no effect for $N\ge29$.  Tables \ref{vector2} and \ref{vector3} treat the two other cases studied in \cite{EL12a}.  In both cases the Parrondo effect appears for games $A$, $B$, and $C$ as well as for games $A'$, $B$, and $C'$, at least for $N\ge6$, whereas the approximate formulas yield misleading results.

\begin{table}[ht]
\caption{\label{Toral}Analysis of the Parrondo effect for Toral's choice of the probability parameters, $(p_0,p_1,p_2)=(1,4/25,7/10)$.  $\mu_B$, $\mu_C$, and $\mu_{C'}$ are mean profits for games $B$, $C$, and $C'$ with $\gamma=1/2$.  $\hat\mu_B$ and $\hat\mu_{C'}$ are the approximations due to Li et al.\ \cite{L14} of $\mu_B$ and $\mu_{C'}$.  Numbers have been rounded to six significant digits.\medskip}
\catcode`@=\active \def@{\hphantom{0}}
\catcode`#=\active \def#{\hphantom{$-$}}
\begin{center}
\begin{footnotesize}
\begin{tabular}{cccccc}
\hline
\noalign{\smallskip}
 $N$ & $\mu_B$ & $\hat\mu_B$ & $\mu_C$ & $\mu_{C'}$ & $\hat\mu_{C'}$ \\
\noalign{\smallskip}
\hline
\noalign{\smallskip}
@@@3 & $-0.0909091$@@ & $-0.0909091$@@ & $-0.0183774$@& $-0.0766158$@ & $-0.0766158$@ \\
@@@4 &  #0.0799608@@  & $-0.0218156$@@ &  #0.0171357@ &  #0.0156538@  & $-0.0424145$@ \\
@@@5 & $-0.00219465$@ & $-0.0136466$@@ &  #0.00405176 &  #0.00565126  & $-0.0293182$@ \\
@@@6 & $-0.0189247$@@ & $-0.0101518$@@ &  #0.00463310 &  #0.01343312  & $-0.0219988$@ \\
@@@7 &  #0.00350598@  & $-0.00790411$@ &  #0.00482261 &  #0.00680337  & $-0.0172930$@ \\
@@@8 &  #0.000698188  & $-0.00620890$@ &  #0.00479021 &  #0.00678290  & $-0.0140074$@ \\
@@@9 & $-0.00189233$@ & $-0.00484806$@ &  #0.00479036 &  #0.00678314  & $-0.0115820$@ \\
@@10 & $-0.000332809$ & $-0.00372258$@ &  #0.00479099 &  #0.00678338  & $-0.00971779$ \\
@@11 & $-0.000466527$ & $-0.00277480$@ &  #0.00479089 &  #0.00678336  & $-0.00823997$ \\
@@12 & $-0.000676916$ & $-0.00196613$@ &  #0.00479089 &  #0.00678336  & $-0.00703965$ \\
@@13 & $-0.000562901$ & $-0.00126876$@ &  #0.00479089 &  #0.00678336  & $-0.00604536$ \\   
@@14 & $-0.000569340$ & $-0.000661814$ &  #0.00479089 &  #0.00678336  & $-0.00520823$ \\
@@15 & $-0.000586184$ & $-0.000129283$ &  #0.00479089 &  #0.00678336  & $-0.00449372$ \\
@@16 & $-0.000578161$ &  #0.000341368  &  #0.00479089 &  #0.00678336  & $-0.00387672$ \\
@@17 & $-0.000578345$ &  #0.000760068  &  #0.00479089 &  #0.00678336  & $-0.00333856$ \\
@@18 & $-0.000579652$ &  #0.00113478@  &  #0.00479089 &  #0.00678336  & $-0.00286501$ \\
@@19 & $-0.000579095$ &  #0.00147194@  &  #0.00479089 &  #0.00678336  & $-0.00244512$ \\
@@20 &                &  #0.00177683@  &              &               & $-0.00207024$ \\
@100 &                &  #0.00652292@  &              &               & #0.00329682   \\
@500 &                &  #0.00748377@  &              &               & #0.00430074   \\
2500 &                &  #0.00767594@  &              &               & #0.00449892   \\
\noalign{\smallskip}
\hline
\end{tabular}
\end{footnotesize}
\end{center}
\end{table}

\begin{table}[ht]
\caption{\label{vector2}Analysis of the Parrondo effect for a second point $(p_0,p_1,p_2)=(7/10,17/25,0)$ on the boundary of the unit cube. \medskip}
\catcode`@=\active \def@{\hphantom{0}}
\catcode`#=\active \def#{\hphantom{$-$}}
\begin{center}
\begin{footnotesize}
\begin{tabular}{cccccc}
\hline
\noalign{\smallskip}
 $N$ & $\mu_B$ & $\hat\mu_B$ & $\mu_C$ & $\mu_{C'}$ & $\hat\mu_{C'}$\\
\noalign{\smallskip}
\hline
\noalign{\smallskip}
@@@3 &   #0.0710383@  & #0.0710383 & #0.0297791@ & #0.0525560@@ & #0.0525560 \\
@@@4 &  $-0.0425713$@ & #0.0485411 & #0.00241457 & #0.000952648 & #0.0363651 \\
@@@5 &   #0.00257895  & #0.0398300 & #0.00818232 & #0.00765099@ & #0.0295117 \\
@@@6 &  $-0.0102930$@ & #0.0349801 & #0.00721881 & #0.0136825@@ & #0.0256872 \\
@@@7 &  $-0.00722622$ & #0.0318731 & #0.00736816 & #0.00691714@ & #0.0232447 \\
@@@8 &  $-0.00808338$ & #0.0297097 & #0.00734464 & #0.00691038@ & #0.0215492 \\
@@@9 &  $-0.00784318$ & #0.0281159 & #0.00734835 & #0.00691100@ & #0.0203035 \\
@@10 &  $-0.00790952$ & #0.0268928 & #0.00734776 & #0.00691094@ & #0.0193494 \\
@@11 &  $-0.00789119$ & #0.0259243 & #0.00734786 & #0.00691095@ & #0.0185954 \\
@@12 &  $-0.00789624$ & #0.0251385 & #0.00734784 & #0.00691095@ & #0.0179843 \\
@@13 &  $-0.00789485$ & #0.0244881 & #0.00734784 & #0.00691095@ & #0.0174792 \\
@@14 &  $-0.00789523$ & #0.0239408 & #0.00734784 & #0.00691095@ & #0.0170546 \\
@@15 &  $-0.00789513$ & #0.0234740 & #0.00734784 & #0.00691095@ & #0.0166927 \\
@@16 &  $-0.00789516$ & #0.0230711 & #0.00734784 & #0.00691095@ & #0.0163806 \\
@@17 &  $-0.00789515$ & #0.0227198 & #0.00734784 & #0.00691095@ & #0.0161086 \\
@@18 &  $-0.00789515$ & #0.0224108 & #0.00734784 & #0.00691095@ & #0.0158696 \\
@@19 &  $-0.00789515$ & #0.0221369 & #0.00734784 & #0.00691095@ & #0.0156577 \\
@@20 &                & #0.0218925 &             &              & #0.0154688 \\
@100 &                & #0.0183996 &             &              & #0.0127773   \\
@500 &                & #0.0177477 &             &              & #0.0122767   \\
2500 &                & #0.0176191 &             &              & #0.0121780   \\
\noalign{\smallskip}
\hline
\end{tabular}
\end{footnotesize}
\end{center}
\end{table}

\begin{table}[ht]
\caption{\label{vector3}Analysis of the Parrondo effect for a point $(p_0,p_1,p_2)=(1/10,3/5,3/4)$ in the interior of the unit cube.\medskip}
\catcode`@=\active \def@{\hphantom{0}}
\catcode`#=\active \def#{\hphantom{$-$}}
\begin{center}
\begin{footnotesize}
\begin{tabular}{cccccc}
\hline
\noalign{\smallskip}
 $N$ & $\mu_B$ & $\hat\mu_B$ & $\mu_C$ & $\mu_{C'}$ & $\hat\mu_{C'}$\\
\noalign{\smallskip}
\hline
\noalign{\smallskip}
@@@3 &  $-0.190476$@@ & $-0.190476$@@@ & $-0.00671141$ & #0.0250737 & #0.0250737 \\
@@@4 &  $-0.0858189$@ & $-0.141636$@@@ &  #0.0108365@  & #0.0175362 & #0.0217807 \\
@@@5 &  $-0.0389980$@ & $-0.106662$@@@ &  #0.0141217@  & #0.0169208 & #0.0202632 \\
@@@6 &  $-0.0183165$@ & $-0.0807523$@@ &  #0.0147166@  & #0.0336654 & #0.0193758 \\
@@@7 &  $-0.00924232$ & $-0.0609171$@@ &  #0.0148223@  & #0.0168224 & #0.0187918 \\
@@@8 &  $-0.00528548$ & $-0.0453218$@@ &  #0.0148408@  & #0.0168213 & #0.0183779 \\
@@@9 &  $-0.00356984$ & $-0.0327906$@@ &  #0.0148441@  & #0.0168212 & #0.0180692 \\
@@10 &  $-0.00282963$ & $-0.0225375$@@ &  #0.0148446@  & #0.0168212 & #0.0178301 \\
@@11 &  $-0.00251155$ & $-0.0140186$@@ &  #0.0148447@  & #0.0168211 & #0.0176394 \\
@@12 &  $-0.00237531$ & $-0.00684626$@ &  #0.0148447@  & #0.0168211 & #0.0174837 \\
@@13 &  $-0.00231709$ & $-0.000737736$ &  #0.0148448@  & #0.0168211 & #0.0173543 \\
@@14 &  $-0.00229226$ &  #0.00451782@  &  #0.0148448@  & #0.0168211 & #0.0172449 \\
@@15 &  $-0.00228169$ &  #0.00908041@  &  #0.0148448@  & #0.0168211 & #0.0171513 \\
@@16 &  $-0.00227719$ &  #0.0130734@@  &  #0.0148448@  & #0.0168211 & #0.0170703 \\
@@17 &  $-0.00227528$ &  #0.0165933@@  &  #0.0148448@  & #0.0168211 & #0.0169994 \\
@@18 &  $-0.00227446$ &  #0.0197165@@  &  #0.0148448@  & #0.0168211 & #0.0169370 \\
@@19 &  $-0.00227412$ &  #0.0225043@@  &  #0.0148448@  & #0.0168211 & #0.0168815 \\
@@20 &                &  #0.0250062@@  &               &            & #0.0168319 \\
@100 &                &  #0.0610154@@  &               &            & #0.0161143   \\
@500 &                &  #0.0674780@@  &               &            & #0.0159784   \\
2500 &                &  #0.0687331@@  &               &            & #0.0159515   \\
\noalign{\smallskip}
\hline
\end{tabular}
\end{footnotesize}
\end{center}
\end{table}

Another distinction is the rate at which the means $\mu_B$, $\mu_C$, and $\mu_{C'}$ converge.  Notice that $\mu_C$ and $\mu_{C'}$ converge very rapidly, having stabilized to six significant digits by $N=11$, 12, or 13 in each case.  $\mu_B$ converges a little more slowly.  By $N=19$, it has stabilized to three significant digits in the first case, six in the second case, and four in the third case.  As a consequence, it is unnecessary to extend these calculations to larger $N$.  If it were possible, the results would be virtually identical.  Contrast that with the situation for the approximate means $\hat\mu_B$ and $\hat\mu_{C'}$ of Li et al.\ \cite{L14}.  As one can see from Tables \ref{Toral}--\ref{vector3}, the convergence is \textit{much} slower.

\begin{figure}[ht]
\centering
\includegraphics[width = 1.45in]{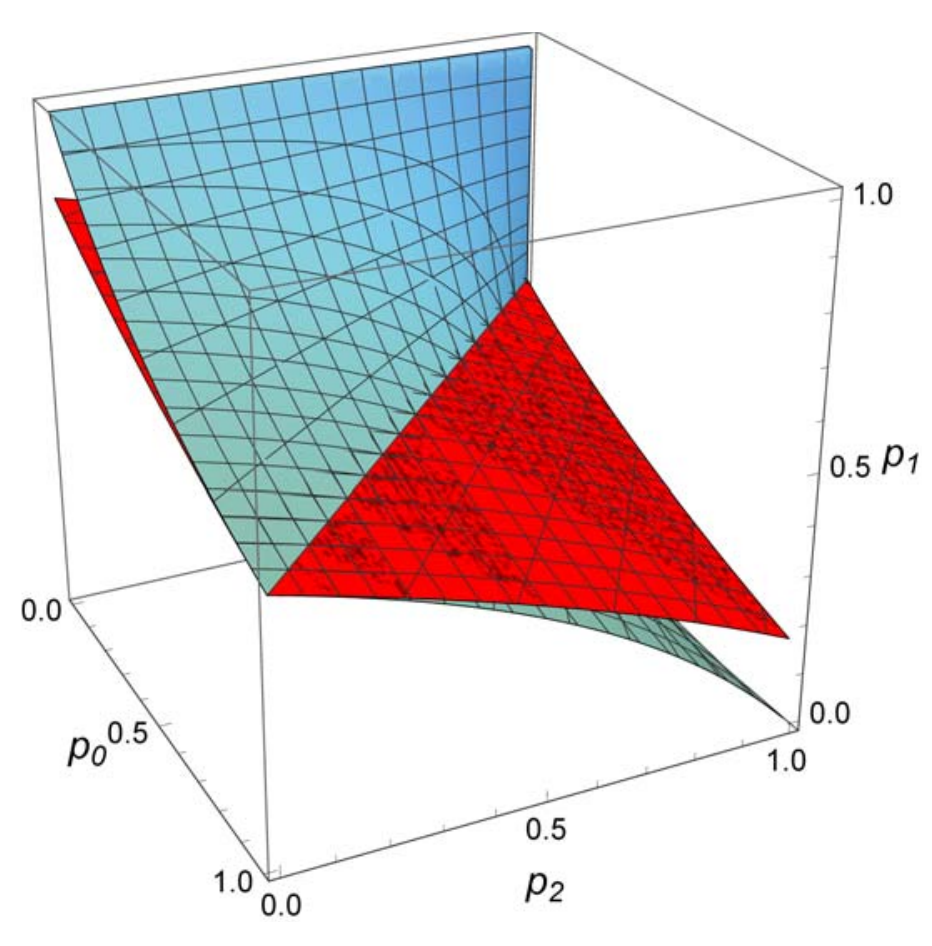}\quad
\includegraphics[width = 1.45in]{N=4}\quad
\includegraphics[width = 1.45in]{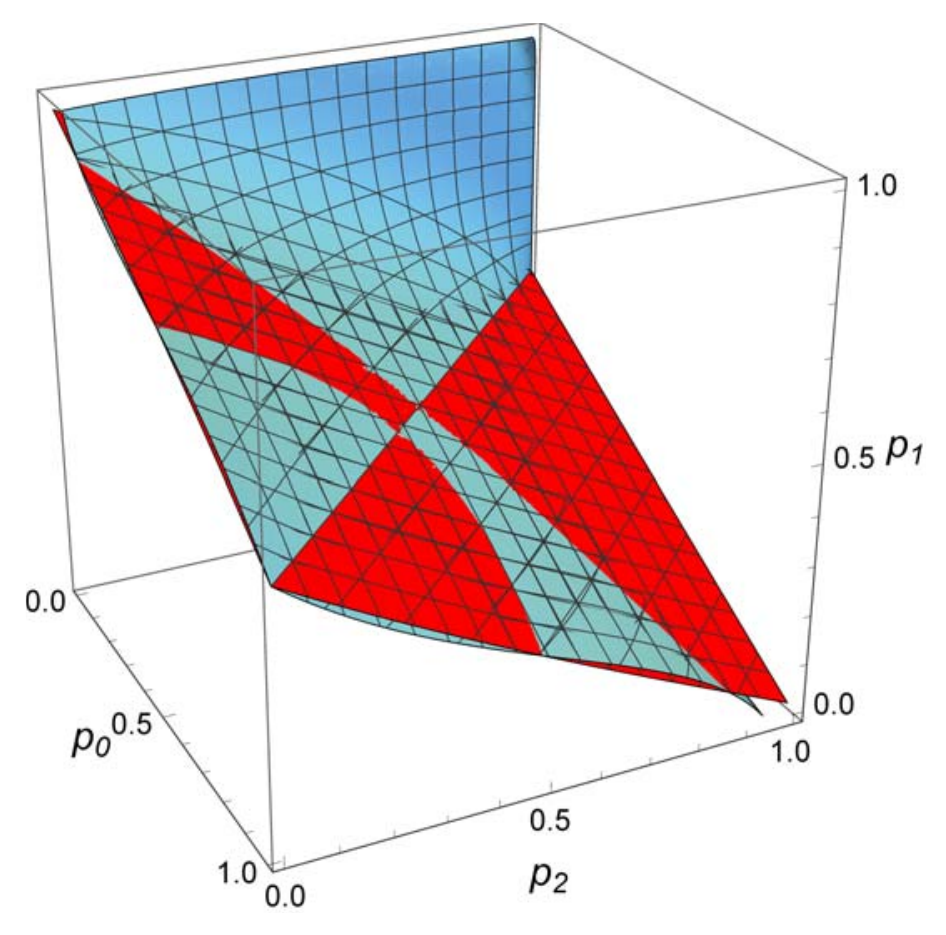}\\
\small $N=3$ \hspace{1.2in} $N=4$  \hspace{1.2in} $N=5$ \\

\bigskip
\includegraphics[width = 1.45in]{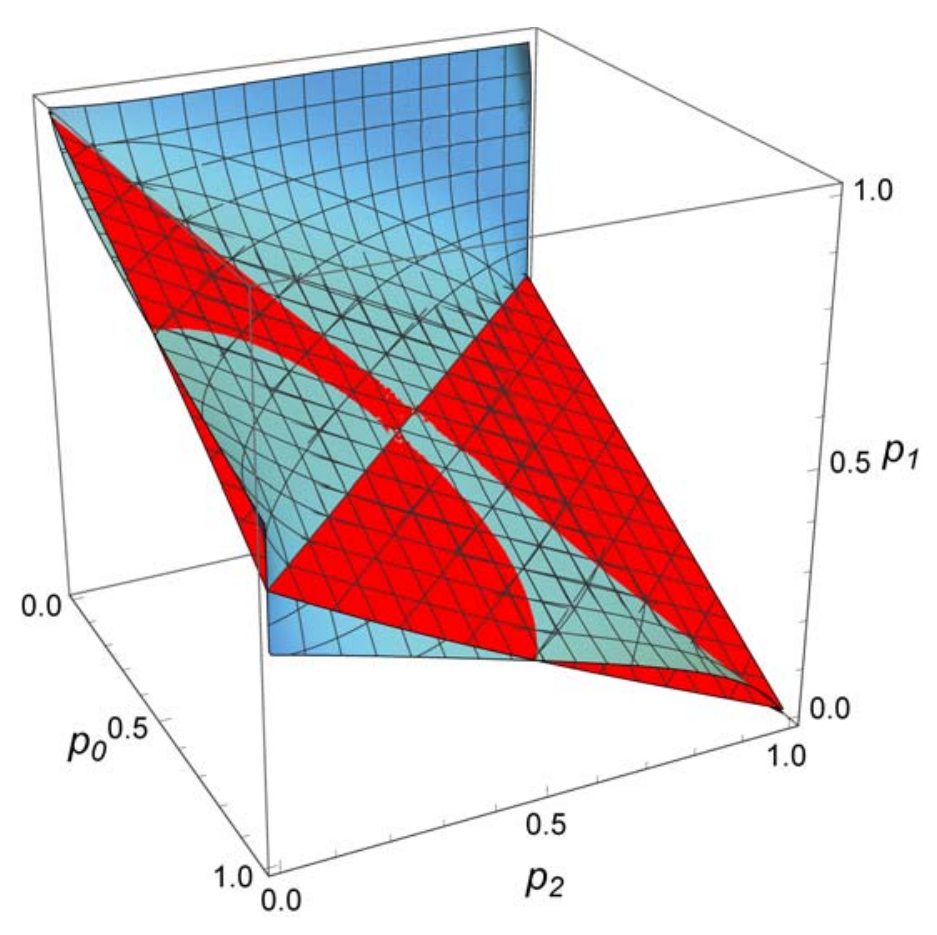}\quad
\includegraphics[width = 1.45in]{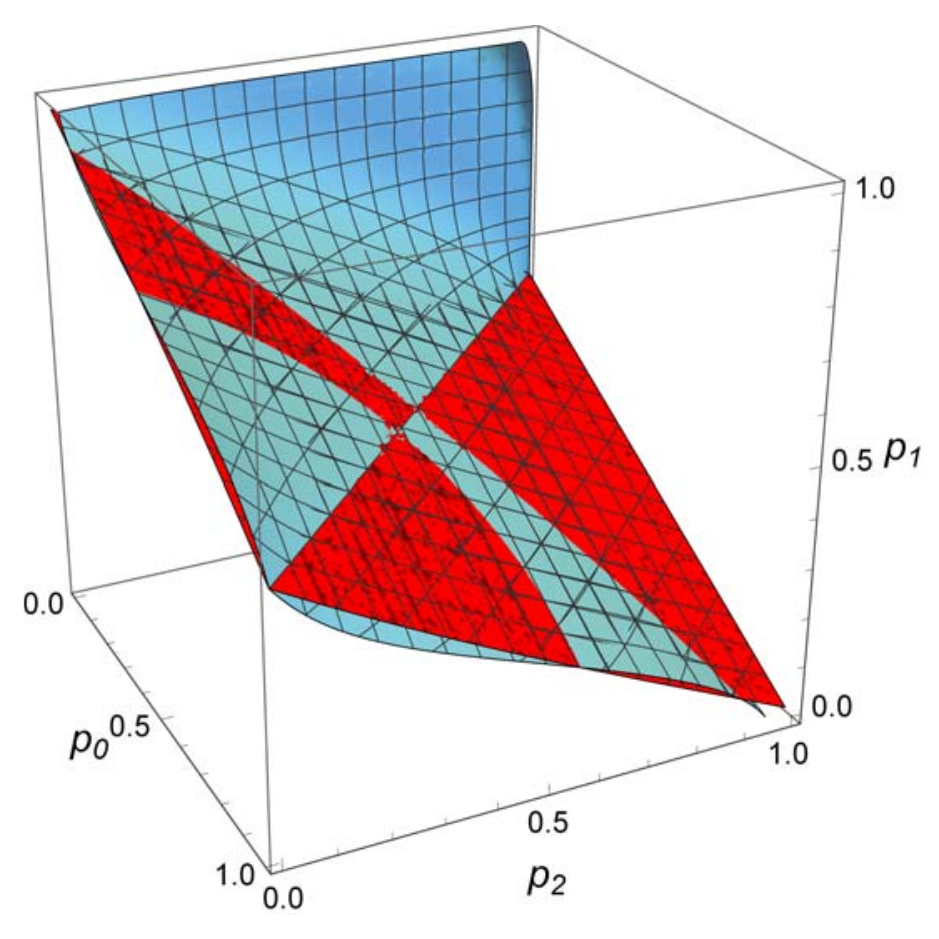}\\
\small $N=6$  \hspace{1.2in} $N=7$  \\

\bigskip
\includegraphics[width = 1.45in]{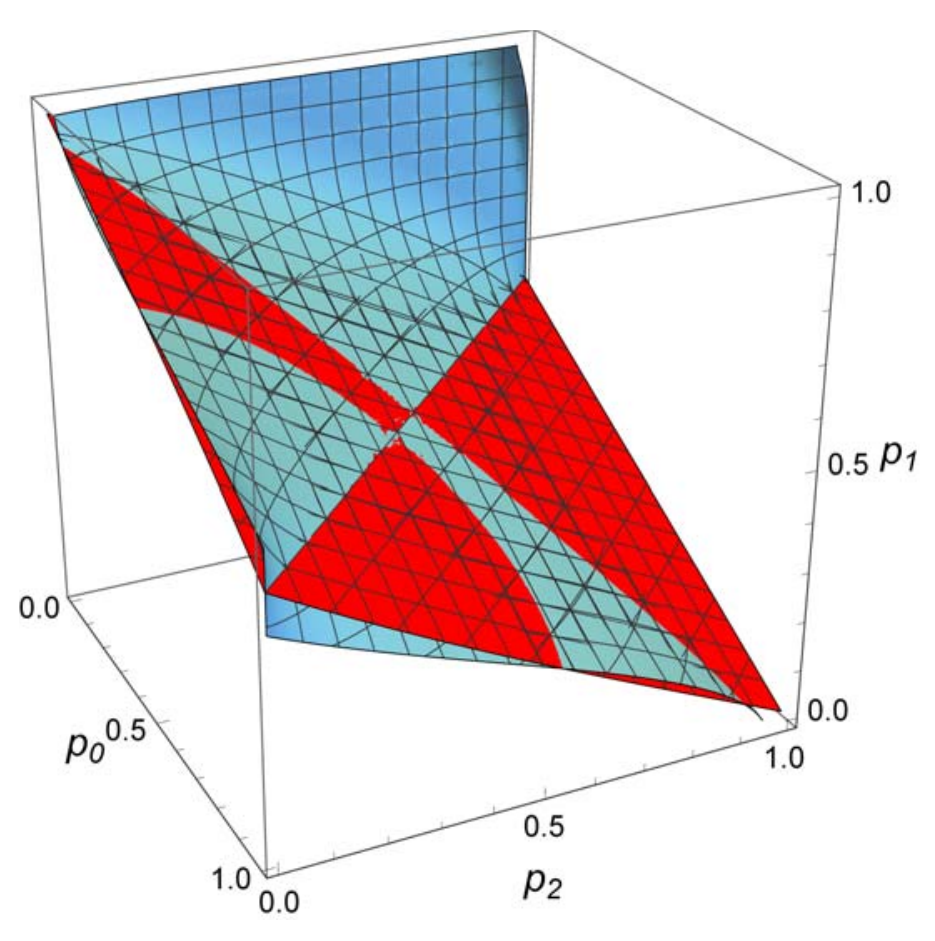}\quad
\includegraphics[width = 1.45in]{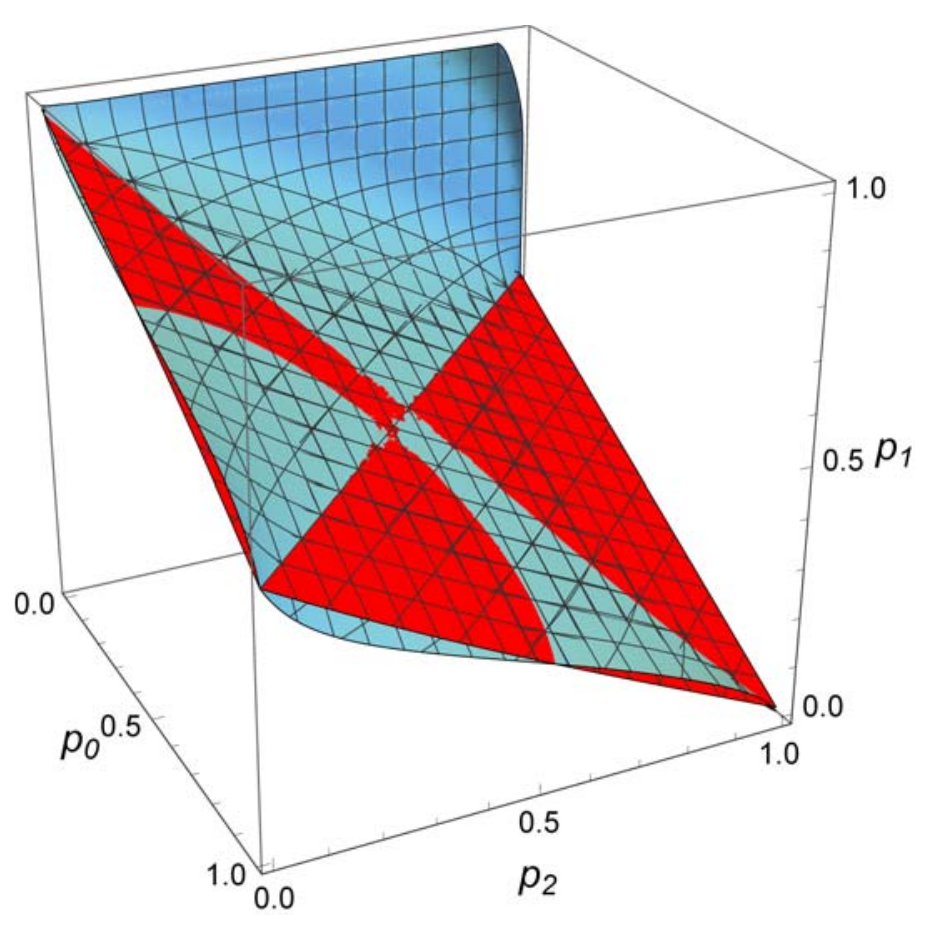}\\
\small $N=8$  \hspace{1.2in} $N=9$ 

\caption{\label{region_EL}For $3\le N\le 9$ and $\gamma=1/2$, the blue surface is the surface $\mu_B=0$, and the red surface is the surface $\mu_{C'}=0$, in the $(p_0,p_2,p_1)$ unit cube.  The Parrondo region is the region on or below the blue surface and above the red surface, while the anti-Parrondo region is the region on or above the blue surface and below the red surface.}
\end{figure}

\begin{figure}[ht]
\centering
\includegraphics[width = 1.45in]{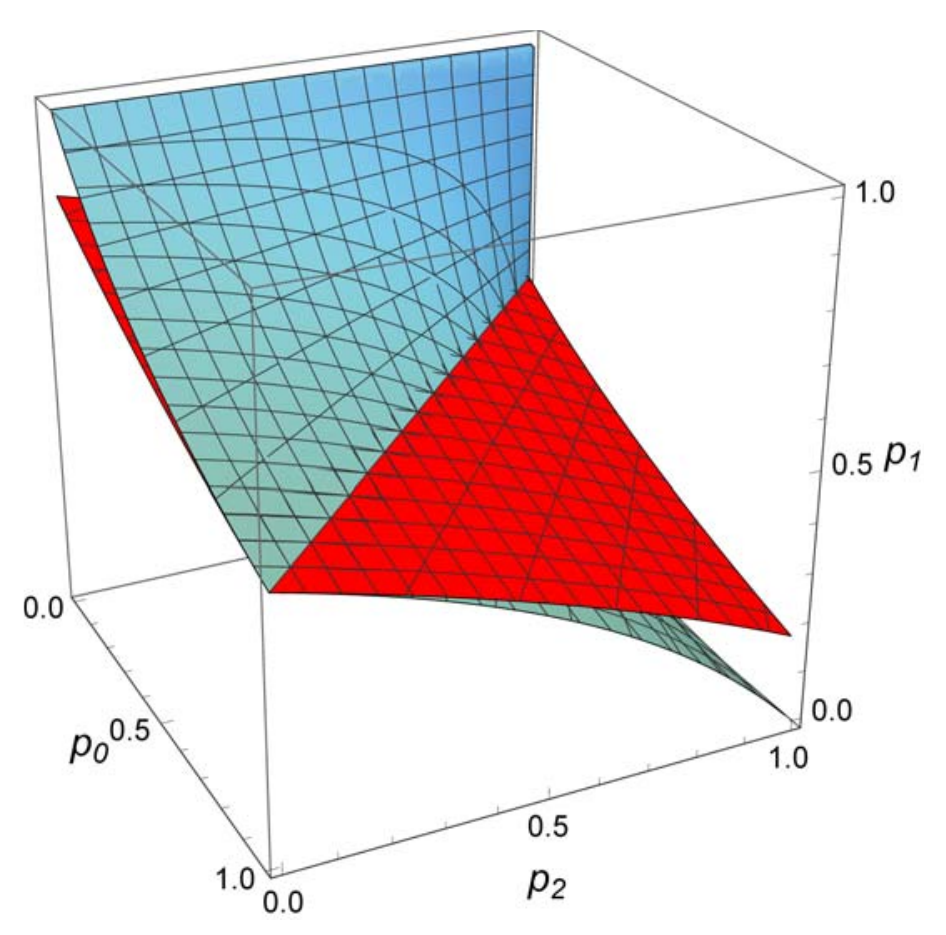}\quad
\includegraphics[width = 1.45in]{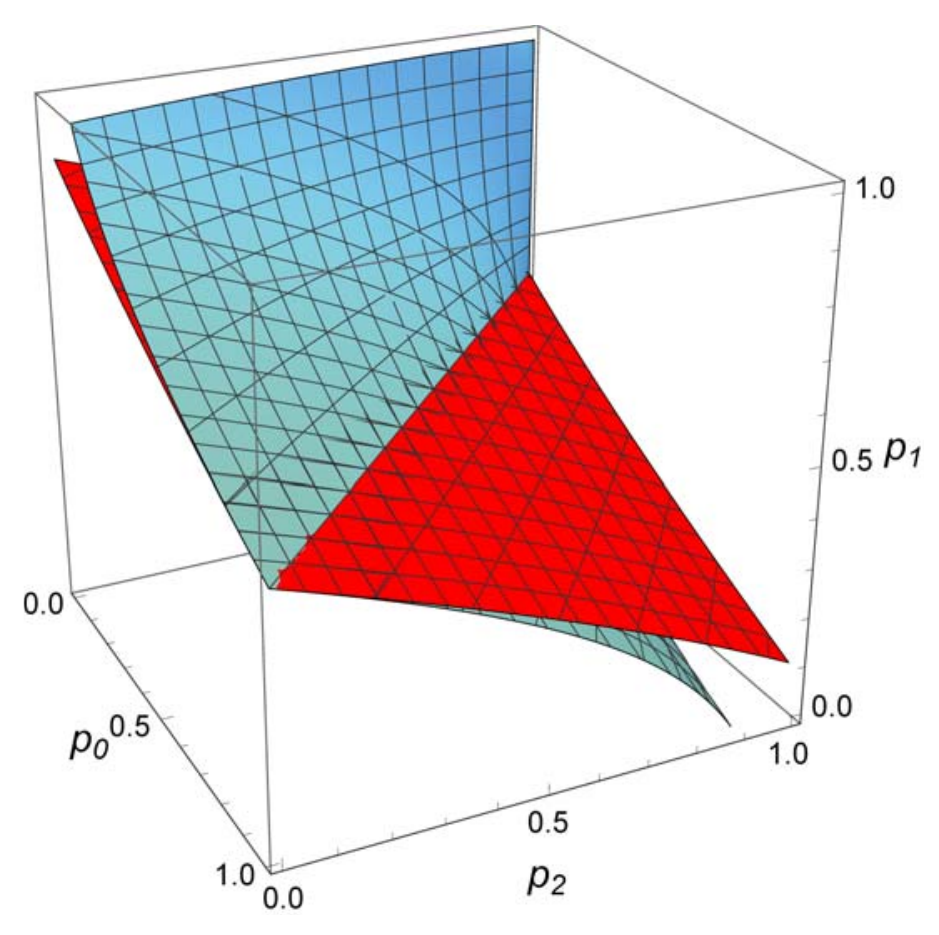}\quad
\includegraphics[width = 1.45in]{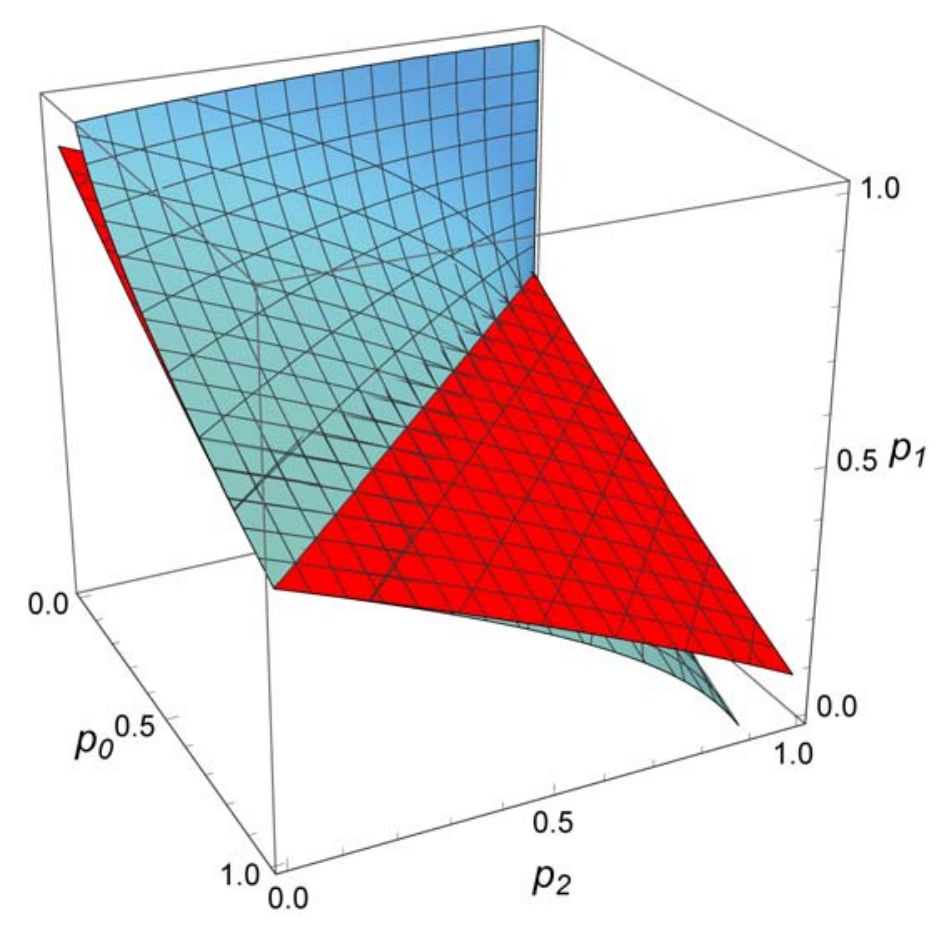}\\
\small $N=3$ \hspace{1.2in} $N=4$  \hspace{1.2in} $N=5$ \\

\bigskip
\includegraphics[width = 1.45in]{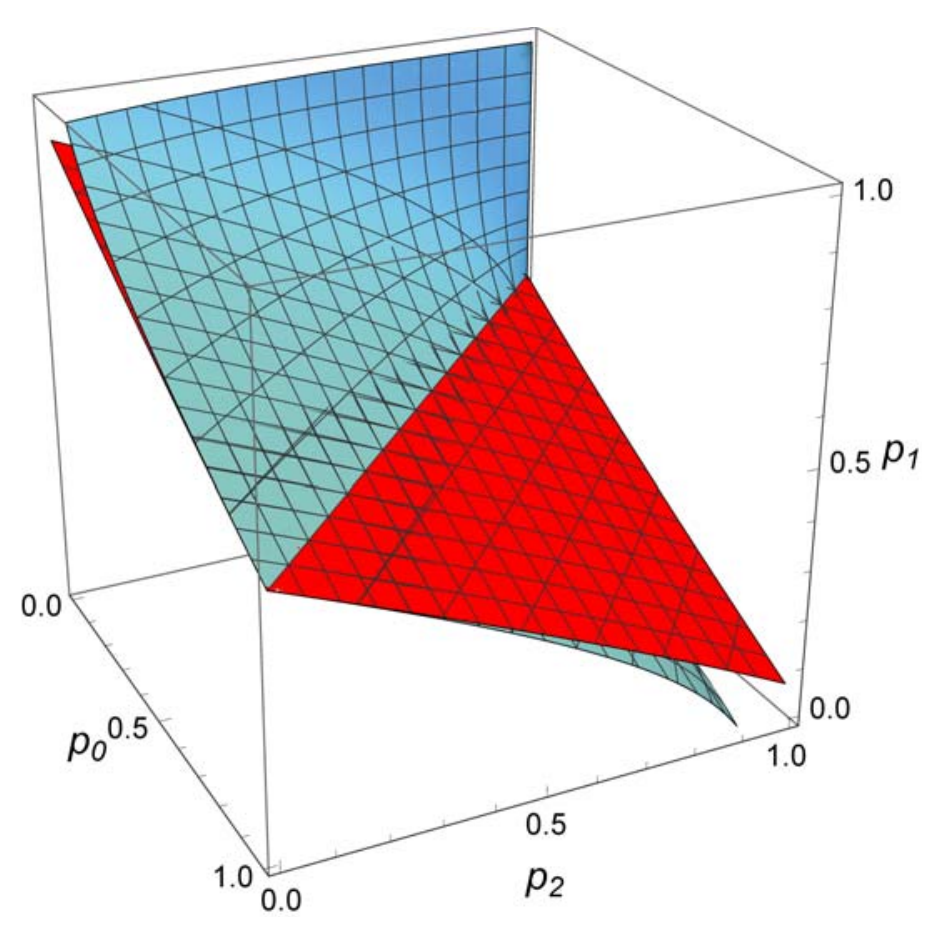}\quad
\includegraphics[width = 1.45in]{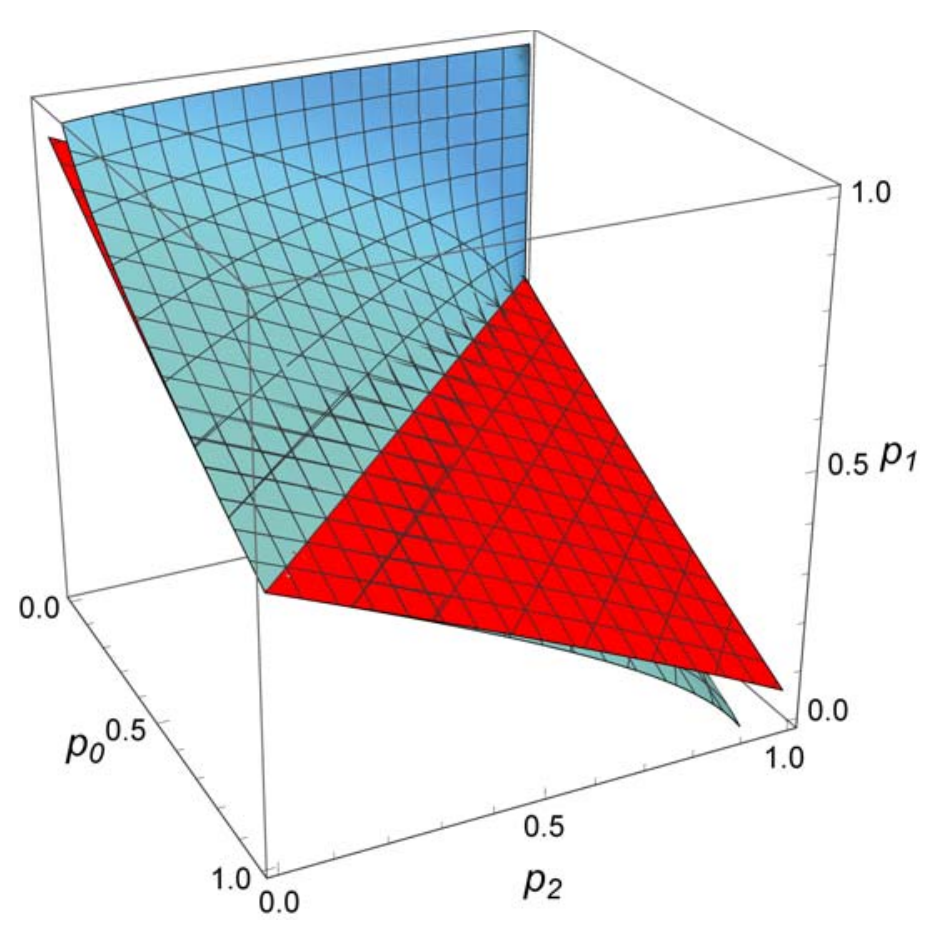}\quad
\includegraphics[width = 1.45in]{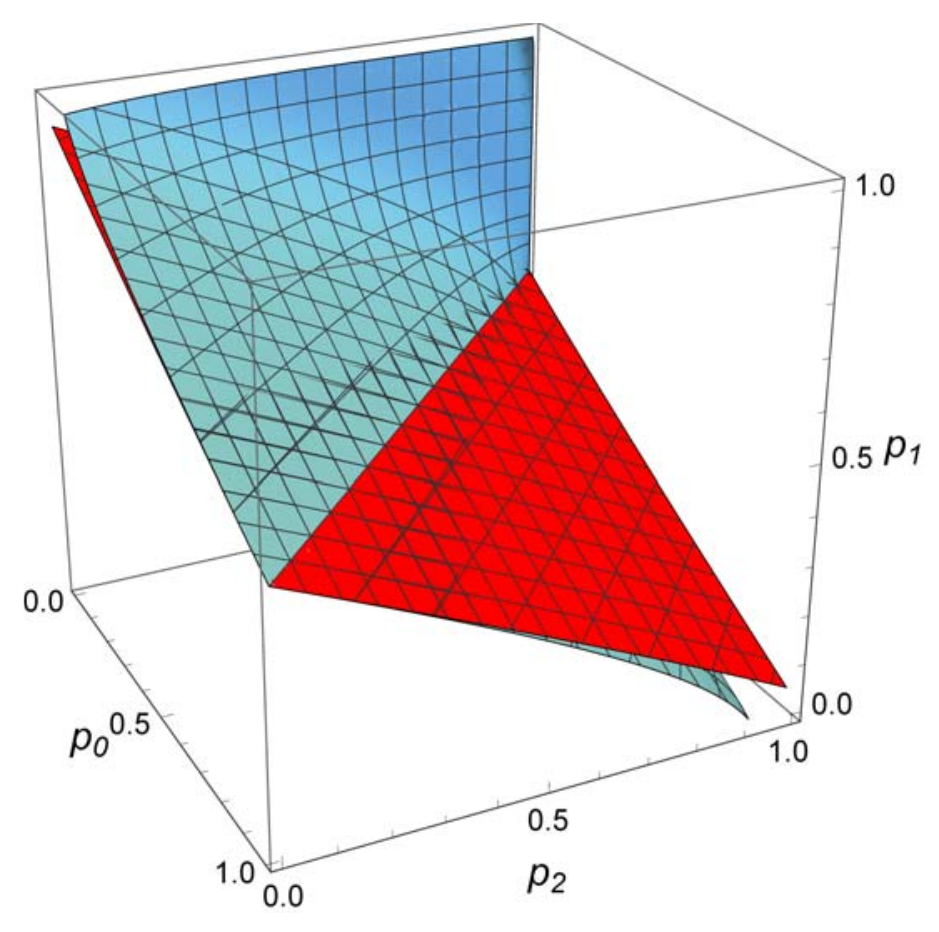}\\
\small $N=6$ \hspace{1.2in} $N=7$  \hspace{1.2in} $N=8$ \\

\bigskip
\includegraphics[width = 1.45in]{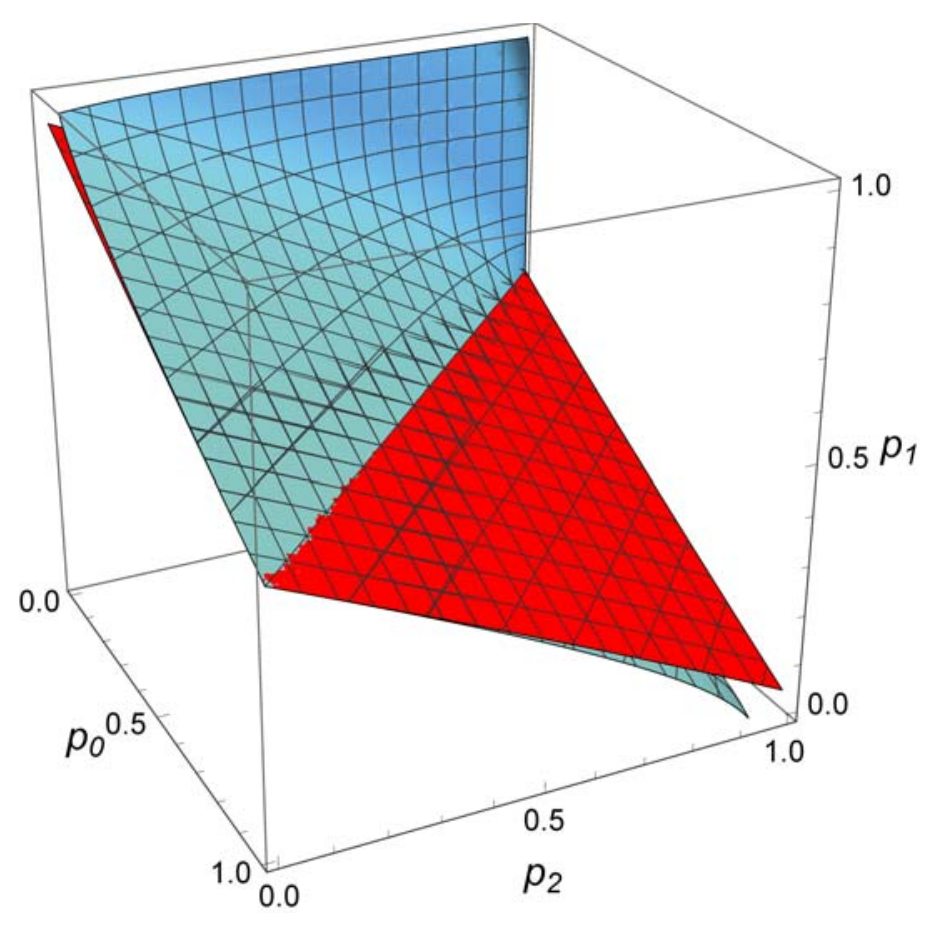}\quad
\includegraphics[width = 1.45in]{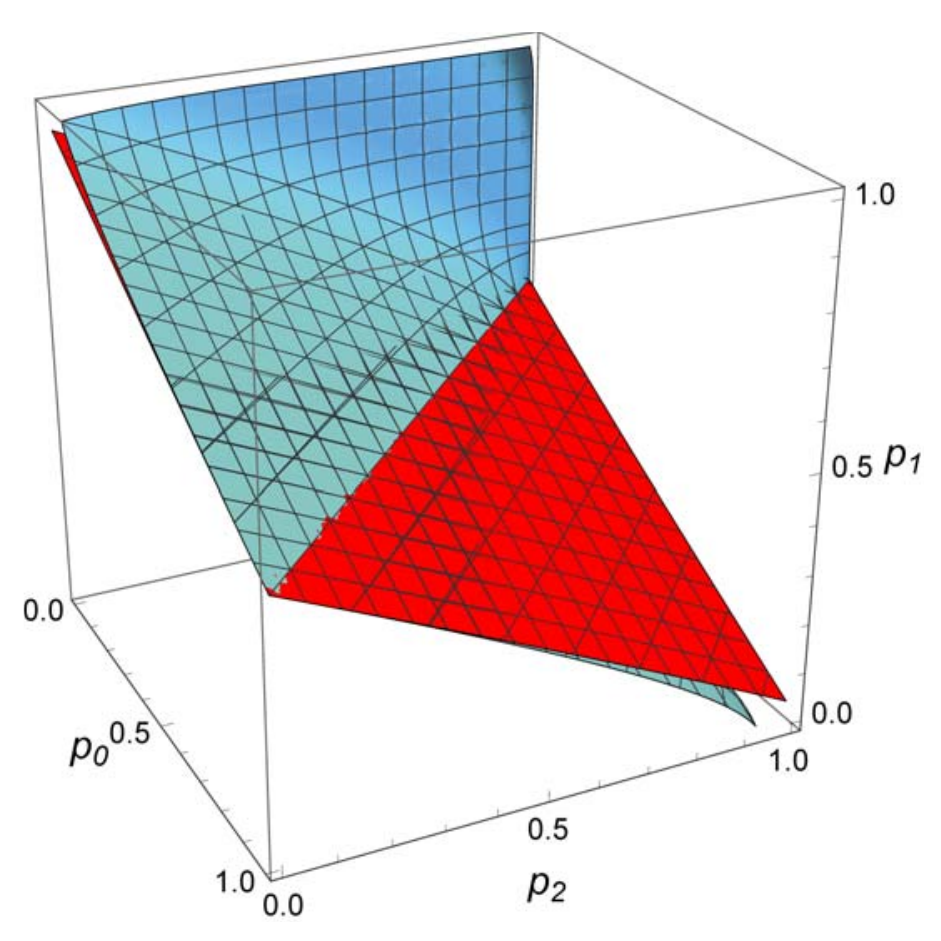}\quad
\includegraphics[width = 1.45in]{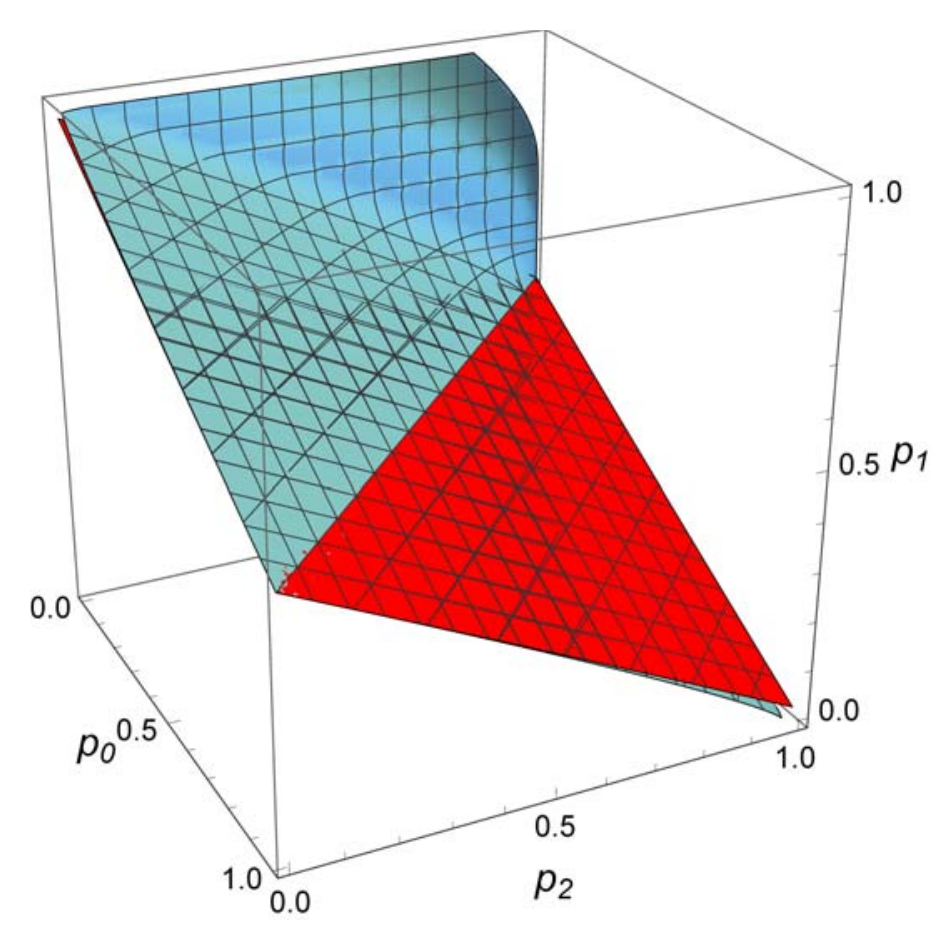}\\
\small $N=9$ \hspace{1.2in} $N=10$  \hspace{1.2in} $N=20$ \\

\caption{\label{region_Li}For $3\le N\le 10$, $N=20$, and $\gamma=1/2$, the blue surface is the surface $\hat\mu_B=0$, and the red surface is the surface $\hat\mu_{C'}=0$, in the $(p_0,p_2,p_1)$ unit cube, where $\hat\mu_B$ and $\hat\mu_{C'}$ are the approximate means of Li et al.\ \cite{L14}. The (approximate) Parrondo region is the region on or below the blue surface and above the red surface, while the (approximate) anti-Parrondo region is the region on or above the blue surface and below the red surface.}
\end{figure}

Next, we use computer graphics to sketch, for games $A'$, $B$, and $C'$, the Parrondo and anti-Parrondo regions when $3\le N\le 9$.  See Figure \ref{region_EL}. We can compare these figures with those of \cite{EL12a}, for games $A$, $B$, and $C$, when $3\le N\le6$.  (The reason we can go further in the present case than we did in \cite{EL12a} is that we found a new method for sketching graphs without the need for explicit formulas.)  The figures for games $A'$, $B$, and $C'$ are distinctively different from those for games $A$, $B$, and $C$.  In both cases, the general shape of the Parrondo and anti-Parrondo regions does not change much, once $N\ge5$.  In Figure \ref{region_Li}, we sketch, for games $A'$, $B$, and $C'$, the approximate Parrondo and anti-Parrondo regions when $3\le N\le 10$ and $N=20$, based on the approximation of Li et al.\ \cite{L14}.  We see clearly that the approximation is poor.

\section{Conclusions}\label{Conclusions}

\begin{itemize}

\item Parrondo games with $N\ge3$ players and one-dimensional spatial dependence were introduced by Toral \cite{T01} and studied analytically by Mihailovi\'c and Rajkovi\'c \cite{MR03} for $N\le12$ and by Ethier and Lee \cite{EL12a} for $N\le19$.  The additional seven cases in the latter work led to the conjecture that the mean profits converge as $N\to\infty$, and this was subsequently proved in \cite{EL13a} under certain conditions.  The reason that larger $N$ could be treated was not because of faster computers but because of a state space reduction method that requires for its justification that a lumpability condition be satisfied by the Markov chains describing the Parrondo games.

\item Xie et al.\ \cite{X11} modified Toral's model by changing the nonspatial game $A$ to one with spatial dependence.  Li et al.\ \cite{L14} studied these Parrondo games using a state space reduction method that reduces the size of the state space from $2^N$ states to $N+1$ states.  However, the Markov chains describing these modified Parrondo games fail to satisfy the lumpability condition, with the result that the reduced Markov chains are inconsistent with the model of Xie et al.  

\item Toral's model as modified by Xie et al.\ \cite{X11} is studied in this paper using the methods introduced in \cite{EL12a}.  The lumpability condition is satisfied and so the reduced Markov chains are consistent with the model of Xie at al.  This allows us to compare our exact results for this model with the approximate ones of Li et al.\ \cite{L14}.  We find that their approximation is poor and that their results are misleading.  This can be seen from Tables \ref{Toral}--\ref{vector3} or by comparing Figures \ref{region_EL} and \ref{region_Li}.  The same is true for large $N$ as well, as Tables \ref{Toral}--\ref{vector3} suggest.

\item Li et al.\ \cite{L14} justified their approximation on the grounds that exact computations are impossible for large $N$.  While that may be true, we would argue that computations for large $N$ are unnecessary.  The Parrondo region stabilizes rather quickly as $N\to\infty$.  By $N=19$, mean profits for the Parrondo games under consideration have stabilized to three or more significant digits.    Although we cannot extrapolate with absolute mathematical certainty, it seems safe to do so based on the computations that have been done so far.  See Tables \ref{Toral}--\ref{vector3}, for example.

\end{itemize}

\end{document}